\documentclass[12pt,leqno]{article}
\usepackage{amssymb}
\usepackage[mathscr]{eucal}
\usepackage{amsmath,amssymb,latexsym,theorem,bbm}
\newcommand{\EE}{\mathsf{E}}

\newcommand{\NN}{\mathbb{N}}
\newcommand{\PP}{\mathsf{P}}

\newcommand{\RR}{\mathbb{R}}

\newcommand{\cB}{{\mathcal B}}

\newcommand{\cF}{{\mathcal F}}

\newcommand{\cL}{{\mathcal L}}
\newcommand{\cN}{{\mathcal N}}

\newcommand{\dd}{\mathrm{d}}
\newcommand{\ee}{\mathrm{e}}

\DeclareMathOperator*{\argmax}{arg\,max}

\newcommand{\sign}{\operatorname{sign}}
\newcommand{\mean}{\stackrel{L_1}{\longrightarrow}}

\newcommand{\halpha}{\widehat{\alpha}}

\newcommand{\tM}{\widetilde{M}}

\renewcommand{\leq}{\leqslant}
\renewcommand{\geq}{\geqslant}

\newcommand{\stoch}{\stackrel{\PP}{\longrightarrow}}
\newcommand{\distr}{\stackrel{\cL}{\longrightarrow}}
\newcommand{\distre}{\stackrel{\cL}{=}}
\newcommand{\qmean}{\stackrel{L_2}{\longrightarrow}}

\newcommand{\vare}{\varepsilon}
\newcommand{\proofend}{\hfill\mbox{$\Box$}}

\numberwithin{equation}{section}

\theoremstyle{change} \theorembodyfont{\em}
\newtheorem{Lem}{Lemma.}[section]
\newtheorem{Thm}[Lem]{Theorem.}

\newtheorem{Cor}[Lem]{Corollary.}

\theorembodyfont{\rm}
\newtheorem{Rem}[Lem]{Remark.}

\begin{document}

\begin{center}
 {\bfseries\Large Asymptotic behavior of maximum likelihood estimator} \\[2mm]
 {\bfseries\Large for time inhomogeneous diffusion processes}\\[5mm]

 {\sc\large M\'aty\'as $\text{Barczy}^{*,\diamond}$} {\large and}
 {\sc\large Gyula $\text{Pap}^*$}
\end{center}

\vskip0.5cm

* University of Debrecen, Faculty of Informatics, Pf.~12, H--4010 Debrecen, Hungary;
 e--mail: barczy@inf.unideb.hu (M. Barczy), papgy@inf.unideb.hu (G. Pap).

$\diamond$ Corresponding author.

\renewcommand{\thefootnote}{}
\footnote{\textit{Acknowledgements\/}:
 The first author has been supported by the Hungarian Scientific Research Fund under
 Grants No.\ OTKA--F046061/2004 and OTKA T-048544/2005.
The second author has been supported by the Hungarian Scientific Research Fund under
 Grant No.\ OTKA T-048544/2005.}

%\vspace*{-10mm}

\begin{abstract}
First we consider a process \ $(X^{(\alpha)}_t)_{t\in[0,T)}$ \ given by a SDE
 \ $\dd X^{(\alpha)}_t = \alpha b(t) X^{(\alpha)}_t \, \dd t + \sigma(t) \, \dd B_t$,
 \ $t \in [0,T)$, \ with a parameter $\alpha \in \RR$, \ where
 \ $T \in (0,\infty]$ \ and \ $(B_t)_{t\in[0,T)}$ \ is a standard Wiener process.
We study asymptotic behavior of the MLE \ $\halpha_t^{(X^{(\alpha)})}$ \ of
 \ $\alpha$ \ based on the observation \ $(X^{(\alpha)}_s)_{s\in[0,\,t]}$ \ as
 \ $t\uparrow T$.
\ We formulate sufficient conditions under which
 \ $\sqrt{I_{X^{(\alpha)}}(t)} \, \big( \halpha_t^{(X^{(\alpha)})} - \alpha \big)$
 \ converges to the distribution of
 \ $c \int_0^1 W_s \, \dd W_s \, \Big/ \int_0^1(W_s)^2 \, \dd s$, \ where
 \ $I_{X^{(\alpha)}}(t)$ \ denotes the Fisher information for \ $\alpha$
 \ contained in the sample \ $(X^{(\alpha)}_s)_{s\in[0,\,t]}$, \ $(W_s)_{s\in [0,1]}$ \
 is a standard Wiener process, and \ $c=1/\sqrt{2}$ \ or \ $c=-1/\sqrt{2}$.
\ We also weaken the sufficient conditions due to Luschgy
 \cite[Section 4.2]{Lus1} under which
 \ $\sqrt{I_{X^{(\alpha)}}(t)} \, \big( \halpha_t^{(X^{(\alpha)})} - \alpha \big)$
 \ converges to the Cauchy distribution.
Furthermore, we give sufficient conditions so that the MLE of \ $\alpha$ \ is
 asymptotically normal with some appropriate random normalizing factor.

Next we study a SDE
 \ $\dd Y^{(\alpha)}_t
    = \alpha b(t) a(Y^{(\alpha)}_t) \, \dd t + \sigma(t) \, \dd B_t$,
 \ $t \in [0,T)$, \ with a perturbed drift satisfying
 \ $a(x) = x + O(1 + |x|^\gamma)$ \ with some \ $\gamma \in [0,1)$.
\ We give again sufficient conditions under which
 \ $\sqrt{I_{Y^{(\alpha)}}(t)} \,
    \big( \widehat\alpha_t^{(Y^{(\alpha)})} - \alpha \big)$
 \ converges to the distribution of
 \ $c \int_0^1W_s \, \dd W_s \, \Big/ \int_0^1(W_s)^2 \, \dd s$.

We emphasize that our results are valid in both cases \ $T\in(0,\infty)$ \ and
 \ $T=\infty$, \ and we develope a unified approach to handle these cases.
\end{abstract}

\vskip1cm

\section{Introduction}
\label{Intro}

Statistical estimation of parameters of diffusion processes has been studied
 for a long time.
Feigin \cite{Fei} gave a good historical overview of the very early
 investigations and provided a general asymptotic theory of maximum likelihood
 estimation (MLE) for continuous-time homogeneous diffusion processes without
 stationarity assumptions and without to resorting to the use of stopping
 times.
Feigin \cite{Fei} also demonstrated the role of martingale limit theory in the
 theory of statistical inference for stochastic processes.
Since then the problem of estimating the drift parameter based on continuous
 observations of time homogeneous diffusions has been extensively studied, see,
 e.g., the books of Liptser and Shiryaev \cite{LipShiI}, \cite{LipShiII} and Kutoyants \cite{Kut}.
For time inhomogeneous diffusions, we can address the books of
 Basawa and Prakasa Rao \cite{BasRao}, Kutoyants \cite{Kut3} and Bishwal \cite{Bis},
 and the research paper of Mishra and Prakasa Rao \cite{MisRao}.

Let \ $T \in (0,\infty]$ \ be fixed.
Let us consider a time inhomogeneous diffusion process
 \ $(Y_t^{(\alpha)})_{t \in [0,T)}$ \ given by the stochastic differential equation
 (SDE)
  \begin{equation}\label{perturbed_SDE}
  \begin{cases}
   \dd Y_t^{(\alpha)}
   = \alpha b(t) a(Y_t^{(\alpha)}) \, \dd t
     + \sigma(t) \, \dd B_t, \qquad t \in [0,T) , \\
   \phantom{\dd} Y_0^{(\alpha)} = 0 ,
  \end{cases}
 \end{equation}
 where \ $b:[0,T)\to\RR$, \ $a:\RR\to\RR$ \ and \ $\sigma:[0,T)\to(0,\infty)$
 \ are known Borel-measurable functions, \ $(B_t)_{t\in[0,T)}$ \ is a standard
 Wiener process, and \ $\alpha\in\RR$ \ is an unknown parameter.

One can obtain sufficient conditions for asymptotic normality in case
 \ $T = \infty$ \ from the general Theorem 5.1 in Chapter 9 due to Basawa and Prakasa Rao
 \cite{BasRao}, \ namely, if \ $\alpha$, \ $b$, \ $a$ \ and \ $\sigma$ \ are
 such that there exists a unique strong solution of the SDE
 \eqref{perturbed_SDE} \ and
 \begin{align}\label{SEGED_BASAWA_RAO}
  \frac{1}{t}
  \int_0^t \frac{ b(s)^2 a(Y^{(\alpha)}_s)^2 }{\sigma(s)^2} \, \dd s
  \stoch K_\alpha
  \qquad \text{as \ $t\to\infty$,}
 \end{align}
 with some \ $K_\alpha \in(0,\infty)$, \ where \ $\stoch$ \ denotes convergence
 in probability, then the MLE \ $\halpha_t^{(Y^{(\alpha)})}$ \ of \ $\alpha$
 \ based on the observation \ $(Y^{(\alpha)}_s)_{s\in[0,\,t]}$ \ is weakly consistent,
 and \ $\sqrt{t} \, \big( \halpha_t^{(Y^{(\alpha)})} - \alpha \big)$
 \ converges in distribution to the normal distribution with mean \ $0$ \ and
 with variance \ $K_\alpha^{-1}$ \ as \ $t\uparrow \infty$.
\ We note that this theorem of Basawa and Prakasa Rao \cite{BasRao} is valid for
 multidimensional diffusion processes and the drift and diffusion coefficients can
 have a more general form.
It is not easy to check condition \eqref{SEGED_BASAWA_RAO}, and hence, as a
 general task, it is desirable to describe the asymptotic behavior of the MLE
 of \ $\alpha$ \ (considering more general normalizing factor than \ $\sqrt{t}$)
 \ by giving simpler sufficient conditions.

In the first part of the present paper we investigate the SDE
 \eqref{perturbed_SDE} with \ $a(x)=x$, \ $x\in\RR$, \ namely,
 \begin{equation}\label{special_SDE}
  \begin{cases}
   \dd X_t^{(\alpha)}
   = \alpha b(t) X_t^{(\alpha)} \, \dd t
     + \sigma(t) \, \dd B_t, \qquad t \in [0,T) , \\
   \phantom{\dd} X_0^{(\alpha)} = 0 ,
  \end{cases}
 \end{equation}
 which is a special case of Hull-White (or extended Vasicek) model, see, e.g.,
 Bishwal \cite[page 3]{Bis}.
As one of our main results, we give sufficient conditions under which the MLE
 of \ $\alpha$
 \ normalized by Fisher information converges to the distribution of
 \ $c \int_0^1 W_s \, \dd W_s \, \Big/ \int_0^1(W_s)^2 \, \dd s$, \ where
 \ $(W_s)_{s\in[0,1]}$ \ is a standard Wiener process, and \ $c=1/\sqrt{2}$ \ or
 \ $c=-1/\sqrt{2}$, \ see Theorem \ref{sing_special}.
In the special case \ $T=\infty$ \ and \ $\sigma\equiv1$ \ Luschgy
 \cite[Section 4.2]{Lus1} gave conditions for the MLE of \ $\alpha$
 \ normalized by Fisher information to converge to
 a normal or to a Cauchy distribution.
In case of Cauchy limit distribution, we weaken and generalize conditions of
 Luschgy, see Theorem \ref{PRO_CAUCHY} and Remark \ref{REM_LUSCHGY1}.
Moreover, one can easily formulate conditions for asymptotic normality
 generalizing Luschgy's conditions, see Theorem \ref{PRO_NORMAL}.
(We do not know whether any other limit distribution can appear.)
We also prove that, under the conditions of Theorem \ref{PRO_CAUCHY}
 or Theorem \ref{PRO_NORMAL}, the MLE of \ $\alpha$ \ is asymptotically normal
 with an appropriate {\sl random} normalizing factor,
 see Corollaries \ref{COR_CAUCHY} and \ref{COR_NORMAL}.
Furthermore, we prove strong consistency of the MLE of \ $\alpha$, \ see
 Theorem \ref{PRO_strongly_consistent2}.

The above results are generalizations of the case of an Ornstein-Uhlenbeck
 process, when \ $T=\infty$, \ $b\equiv1$, \ $a(x)=x$, \ $x\in\RR$, \ and
 \ $\sigma\equiv1$.
In this special case if \ $\alpha<0$, \ then the MLE of \ $\alpha$ \ is
 asymptotically normal.
This fact is known for a long time, see, e.g., Example 1.35 in Kutoyants
 \cite{Kut}, (1.3) in Dietz and Kutoyants \cite{DieKut}, page 189 in Basawa and
 Prakasa Rao \cite{BasRao} or Example 2.1 in Gushchin \cite{Gus}.
If \ $\alpha>0$, \ then the MLE of \ $\alpha$ \ is asymptotically Cauchy.
This result is also known for a while, see, e.g., Basawa and Scott
 \cite{BasSco}, Kutoyants \cite{Kut2}, Theorem 5.1 in Dietz and Kutoyants
 \cite{DieKut} or Example 2.1 in Gushchin \cite{Gus}.
If \ $\alpha=0$, \ then
 \[
   t\widehat\alpha_t^{(X^{(0)})}
    \distre \frac{\int_0^1 W_s \, \dd W_s}{\int_0^1 (W_s)^2 \, \dd s},
      \qquad t\in(0,\infty) ,
 \]
 where \ $\distre$ \ denotes equality in distribution, and hence we have not only
 a limit theorem but the appropriately normalized MLE of \ $\alpha$ \ has the
 same distribution for all \ $t \in (0,\infty)$.
\ This has also been known for a long time, see, e.g., (1.4) in Dietz and
 Kutoyants \cite{DieKut}, page 189 in Basawa and Prakasa Rao \cite{BasRao} or
 Example 2.1 in Gushchin \cite{Gus}.
We also note that this distribution is the same as the limit distribution of
 the Dickey-Fuller statistics, see, e.g., the Ph.D. Thesis of Bobkoski
 \cite{Bob}, or (7.14) and Theorem 9.5.1 in Tanaka \cite{Tan}.
The strong consistency of the MLE of \ $\alpha$ \ has also been known for a
 long time, see, e.g., Theorem 17.4 in Liptser and Shiryaev \cite{LipShiII}.

In the second part of the present paper we investigate the SDE
 \eqref{perturbed_SDE}
 with \ $a(x)=x+r(x)$, \ $x\in\RR$ \ and a known Lipschitz function \ $r$
 \ satisfying \ $r(x) = O(1 + |x|^\gamma)$ \ with some \ $\gamma \in [0,1)$,
 \ which can be considered as a perturbation of the SDE \eqref{special_SDE}.
We give sufficient conditions under which the MLE of \ $\alpha$ \ normalized
 by Fisher information converges to the distribution of
 \ $c \int_0^1 W_s \, \dd W_s \, \Big/ \int_0^1(W_s)^2 \, \dd s$, \ where
 \ $c=1/\sqrt{2}$ \ or \ $c=-1/\sqrt{2}$, \ see Theorem \ref{sing_perturb}.
Our proof is based on a generalization of Gr\"onwall's inequality
 (see, Lemma \ref{Gronwall}).
Note that Dietz and Kutoyants \cite{DieKut} investigated the asymptotic
 properties of the MLE \ $\widehat\alpha_t^{(Y^{(\alpha)})}$ \ of \ $\alpha$ \
 in the special case \ $T=\infty$, \ $\alpha>0$, \ $b(t)=c$, \ $t\geq0$, \ with
 some \ $c>0$, \ and \ $\sigma\equiv1$.
\ They showed
 \[
   \sqrt{\frac{c}{2\alpha}}
   \ee^{\alpha c t}
   \Big(\widehat\alpha_t^{(Y^{(\alpha)})}-\alpha\Big)
   \distr\frac{\xi}{\eta^{(\alpha)}}
   \qquad \text{as \ $t\to \infty$,}
 \]
 where \ $\distr$ \ denotes convergence in distribution,
 \[
   \eta^{(\alpha)}
   := \int_0^\infty \ee^{-\alpha c s} \, \dd B_s
     + \alpha c \int_0^\infty \ee^{-\alpha c s} r(Y_s^{(\alpha)}) \, \dd s ,
 \]
 and \ $\xi$ \ is a standard normally distributed random variable independent
 of \ $\eta^{(\alpha)}$, \ provided that \ $\PP(\eta^{(\alpha)}=0)=0$.
\ Dietz and Kutoyants \cite[Theorem 4.1]{DieKut} also showed that
 \ $\widehat\alpha_t^{(Y^{(\alpha)})}$ \ is strongly consistent provided that
 \ $\PP(\eta^{(\alpha)}=0)=0$.

We emphasize that our results are valid in both cases \ $T\in(0,\infty)$ \ and
 \ $T=\infty$, \ and we develope a unified approach to handle these cases.

\vskip1cm

\section{A special time inhomogeneous SDE}\label{special_SDE_section}

Let \ $T \in (0, \infty]$ \ be fixed. Let \ $b: [0,T) \to \RR$ \ and
 \ $\sigma: [0,T) \to \RR$ \ be continuous functions.
 Suppose that \ $\sigma(t) > 0$ \ for all \ $t\in[0,T)$, \ and there exists
 \ $t_0 \in (0,T)$ \ such that \ $b(t)\ne 0$ \ for all \ $t\in[t_0,T).$ \
 For all \ $\alpha \in \RR$, \ consider the SDE \eqref{special_SDE}.
Note that the drift and diffusion coefficients of the SDE \eqref{special_SDE}
 satisfy the local Lipschitz condition and the linear growth condition
 (see, e.g., Jacod and Shiryaev \cite[Theorem 2.32, Chapter III]{JacShi}).
By Jacod and Shiryaev \cite[Theorem 2.32, Chapter III]{JacShi}, the SDE
 \eqref{special_SDE} has a unique strong solution
 \begin{align}\label{SolX}
   X^{(\alpha)}_t
   = \int_0^t
      \sigma(s) \exp\left\{ \alpha \int_s^t b(u) \, \dd u \right\} \dd B_s ,
   \qquad t \in [0,T) ,
 \end{align}
 defined on a filtered probability space
 \ $\big(\Omega, \cF, (\cF_t)_{t\in[0,T)}, \PP\big)$ \ constructed by the help of
 the standard Wiener process \ $B$, \ see, e.g., Karatzas and Shreve \cite[page 285]{KarShr}.
This filtered probability space satisfies the so called usual conditions, i.e.,
 \ $(\Omega,\cF,\PP)$ \ is complete, the filtration \ $(\cF_t)_{t\in[0,T)}$ \ is right-continuous,
 \ $\cF_0$ \ contains all the $\PP$-null sets in \ $\cF$ \ and \ $\cF=\cF_{T-},$ \ where
 \ $\cF_{T-}:=\sigma\left(\bigcup_{t\in[0,T)}\cF_t\right)$.
\ Note that \ $(X_t^{(\alpha)})_{t \in [0,T)}$ \ has continuous sample paths by the
 definition of strong solution, see, e.g., Jacod and Shiryaev
 \cite[Definition 2.24, Chapter III]{JacShi}.
For all \ $\alpha \in \RR$ \ and \ $t \in (0,T)$, \ let \ $\PP_{X^{(\alpha)},\,t}$
 \ denote the distribution of the process \ $(X_s^{(\alpha)})_{s \in [0,\,t]}$ \ on
 \ $\big(C([0,t]),\cB(C([0,t]))\big)$, \ where \ $C([0,t])$ \ and \ $\cB(C([0,t]))$ \ denote
 the set of all continuous real valued functions defined on \ $[0,t]$ \ and the Borel
 $\sigma$--field on \ $C([0,t])$, \ respectively.
The measures \ $\PP_{X^{(\alpha)},\,t}$ \ and \ $\PP_{X^{(0)},\,t}$ \ are equivalent and
 \begin{align*} %\label{RN}
   \frac{\dd \PP_{X^{(\alpha)},\,t}}{\dd \PP_{X^{(0)},\,t}}
    \left(X^{(\alpha)}\big\vert_{[0,t]}\right)
   = \exp \left\{ \alpha
                  \int_0^t
                   \frac{b(s) \, X_s^{(\alpha)}}{\sigma(s)^2} \, \dd X_s^{(\alpha)}
                  -\frac{\alpha^2}{2}
                   \int_0^t
                    \frac{b(s)^2 \, \big( X_s^{(\alpha)} \big)^2}{\sigma(s)^2}
                     \, \dd s \right\} ,
 \end{align*}
 see Liptser and Shiryaev \cite[Theorem 7.20]{LipShiI}.

For all \ $t\in(0,T)$, \ the maximum likelihood estimator
 \ $\halpha_t^{(X^{(\alpha)})}$ \ of the parameter \ $\alpha$ \ based on the
 observation \ $(X_s^{(\alpha)})_{s \in [0,\,t]}$ \ is defined by
 \[
   \halpha_t^{(X^{(\alpha)})}
    :=\argmax_{\alpha\in\RR}
       \ln\left(\frac{\dd \PP_{X^{(\alpha)},\,t}}{\dd \PP_{X^{(0)},\,t}}
                 \left(X^{(\alpha)}\big\vert_{[0,t]}\right)\right).
 \]

The following lemma guarantees the existence of a unique MLE of \ $\alpha$.

\begin{Lem}\label{LEMMA_SPEC_FELTETEL}
For all \ $\alpha\in\RR$ \ and \ $t\in[t_0,T)$, \ we have
 \[
   \PP\left( \int_0^t
               \frac{b(s)^2 \, \big( X_s^{(\alpha)} \big)^2}{\sigma(s)^2} \, \dd s
               > 0 \right)
   =1.
 \]
\end{Lem}

\noindent{\bf Proof.}
Let \ $\alpha\in\RR$ \ be fixed.
On the contrary, let us suppose that there exists some \ $t_1\in[t_0,T)$ \
 such that \ $\PP(A)>0$, \ where
 \[
   A:=\left\{ \omega \in \Omega :
              \int_0^{t_1}
               \frac{b(s)^2 \, \big( X_s^{(\alpha)}(\omega) \big)^2}
                    {\sigma(s)^2} \, \dd s = 0 \right\}.
 \]
Then for all \ $\omega \in A$, \ we have \ $b(s) X_s^{(\alpha)}(\omega) = 0$
 \ for all \ $s \in [0,t_1]$, \ since \ $b$, \ $\sigma$ \ and
  \ $X_.^{(\alpha)}(\omega)$ \ are continuous on \ $[0,T)$.
\ Using the SDE \eqref{special_SDE}, we get
 \[
   X_s^{(\alpha)}(\omega)
   = \left(\int_0^s \sigma(u) \, \dd B_u\right)(\omega) ,
   \qquad s \in [0,t_1] , \qquad \omega \in A ,
 \]
 and hence
 \[
   b(s)\left(\int_0^s \sigma(u) \, \dd B_u\right)(\omega) = 0 ,
   \qquad s \in [0,t_1] , \qquad \omega \in A .
 \]
By \ $b(t_0) \ne 0$, \ we conclude
 \[
   \PP\left( \int_0^{t_0} \sigma(s) \, \dd B_s = 0 \right) > 0 .
 \]
Here \ $\int_0^{t_0} \sigma(s) \, \dd B_s$ \ is a normally distributed random
 variable with mean \ $0$ \ and with variance
 \ $\int_0^{t_0}\sigma(s)^2 \, \dd s > 0$, \ since \ $\sigma(s) > 0$ \ for all
 \ $s \in [0,T)$, \ which leads us to a contradiction.
\proofend

By Lemma \ref{LEMMA_SPEC_FELTETEL}, for all \ $t\in[t_0,T),$ \ there exists a
 unique maximum likelihood estimator \ $\halpha_t^{(X^{(\alpha)})}$ \ of the
 parameter \ $\alpha$ \ based on the observation \ $(X_s^{(\alpha)})_{s \in [0,\,t]}$
 \ given by
 \[
   \halpha_t^{(X^{(\alpha)})}
   = \frac{\int_0^t
            \frac{b(s) \, X_s^{(\alpha)}}{\sigma(s)^2} \, \dd X_s^{(\alpha)}}
          {\int_0^t
            \frac{b(s)^2 \, ( X_s^{(\alpha)} )^2}{\sigma(s)^2} \, \dd s},
   \qquad  t\in[t_0,T).
 \]
To be more precise, by Lemma \ref{LEMMA_SPEC_FELTETEL}, the maximum likelihood
 estimator \ $\halpha_t^{(X^{(\alpha)})}$, \ $t\in[t_0,T)$, \ exists $\PP$-almost
 surely.
Using the SDE \eqref{special_SDE} we obtain
 \begin{align}\label{SEGED20_uj}
   \halpha_t^{(X^{(\alpha)})} - \alpha
   = \frac{\int_0^t
            \frac{b(s) \, X_s^{(\alpha)}}{\sigma(s)} \, \dd B_s}
          {\int_0^t
            \frac{b(s)^2 \, ( X_s^{(\alpha)} )^2}{\sigma(s)^2} \, \dd s},
   \qquad t\in[t_0,T).
 \end{align}
For all \ $t\in(0,T)$, \ the Fisher information for \ $\alpha$ \ contained in
 the observation \ $(X_s^{(\alpha)})_{s\in[0,\,t]}$, \ is defined by
 \[
   I_{X^{(\alpha)}}(t)
   := \EE \left( \frac{\partial}{\partial\alpha}
                 \ln \left(\frac{\dd \PP_{X^{(\alpha)},\,t}}{\dd \PP_{X^{(0)},\,t}}
                  \left( X^{(\alpha)} \big \vert_{[0,\,t]} \right) \right) \right)^2
   = \int_0^t
      \frac{b(s)^2 \, \EE \big( X_s^{(\alpha)} \big)^2}{\sigma(s)^2}  \, \dd s ,
 \]
 where the last equality follows by the SDE \eqref{special_SDE} and Karatzas and Shreve
 \cite[Proposition 3.2.10]{KarShr}.
Note also that, again by Karatzas and Shreve \cite[Proposition 3.2.10]{KarShr},
 \[
   \EE \big( X_s^{(\alpha)} \big)^2
   = \int_0^s
      \sigma(u)^2
      \exp\left\{ 2 \alpha \int_u^s b(v) \, \dd v \right\} \, \dd u ,
   \qquad s \in [0,T) ,
 \]
 and then, by the conditions on \ $b$ \ and \ $\sigma$,
 \ $\EE\big( X_s^{(\alpha)} \big)^2>0$, \ $s\in(0,T)$, \ and
 \ $I_{X^{(\alpha)}} : (0,T) \to [0,\infty)$ \ is an increasing function with
 \ $I_{X^{(\alpha)}}(t) > 0$ \ for all \ $t \in [t_0,T)$.

The aim of the present paragraph is to formulate a theorem (see Theorem \ref{THM_MCLT_RAO2})
 which we will use for studying asymptotic properties of the MLE of \ $\alpha$.
\ First we recall a limit theorem for continuous local martingales.
Theorem 4.1 in van Zanten \cite{Zan}, which is stated for continuous local martingales with time
 interval \ $[0,\infty)$, \ can be applied to continuous local martingales with time interval
 \ $[0,T)$, \ $T\in(0,\infty)$, \ with appropriate modifications of the conditions, as follows.

\begin{Thm}\label{THM_Zanten}
 Let \ $T\in(0,\infty]$ \ be fixed and let
 \ $\big(\Omega, \cF, (\cF_t)_{t\in[0,T)}, \PP\big)$ \ be a filtered probability space
 satisfying the usual conditions.
Let \ $(M_t)_{t \in [0,T)}$ \ be a continuous local martingale with respect to the
 filtration \ $(\cF_t)_{t\in[0,T)}$ \ such that \ $\PP(M_0=0)=1$.
\ Suppose that there exists a function \ $Q:[0,T)\to \RR\setminus\{0\}$ \ such that
 \ $\lim_{t\uparrow T}Q(t)=0$ \ and
 \[
  Q(t)^2\langle M\rangle_t\stoch \eta^2
    \quad \text{as} \quad  t\uparrow T,
 \]
 where \ $\eta$ \ is a random variable defined on \ $\big(\Omega, \cF, \PP\big)$,
 \ and  \ $(\langle M\rangle_t)_{t\in[0,T)}$ \ denotes the
 quadratic variation of \ $M$.
\ Then for each random variable \ $Z$ \ defined on \ $\big(\Omega, \cF, \PP\big)$,
 \ we have
 \[
       (Q(t)M_t,Z)\distr (\eta\xi,Z) \qquad\text{as}\qquad t\uparrow T,
 \]
 where \ $\xi$ \ is a standard normally distributed random variable independent
 of \ $(\eta,Z)$.
\end{Thm}

\noindent To derive a consequence of Theorem \ref{THM_Zanten} we need the following lemma which is
 a multidimensional version of Lemma 3 due to K\'atai and Mogyor\'odi \cite{KatMog}.

\begin{Lem}\label{LEMMA_KATAI_MOGYORODI}
Let \ $T\in(0,\infty]$ \ be fixed. Suppose that \ $(X_t)_{t\in[0,T)}$ \ and
 \ $(Y_t)_{t\in[0,T)}$ \ are stochastic processes on a probability space
 \ $\big(\Omega, \cF, \PP\big)$ \ such that \ $X_t$ \ converges in distribution
 as \ $t\uparrow T$ \ and \ $Y_t\stoch Y$ \ as \ $t\uparrow T$, \ where
 \ $Y$ \ is a random variable defined on \ $\big(\Omega, \cF, \PP\big)$.
\ If \ $g:\RR^2\to\RR^d$ \ is a continuous function (where \ $d\in\NN$), \ then
 \[
  g(X_t,Y_t)-g(X_t,Y) \stoch 0 \qquad \text{as}\quad t\uparrow T.
 \]
\end{Lem}

\noindent{\bf Proof.}
The assertion follows from K\'atai and Mogyor\'odi \cite[Lemma 3]{KatMog} using that
 convergence in probability of a $d$-dimensional stochastic process is equivalent to
 the convergence in probability of all of its coordinates separately (see, e.g.,
 van der Vaart \cite[page 10]{Vaart}).
\proofend

\noindent As a consequence of Theorem \ref{THM_Zanten} and Lemma \ref{LEMMA_KATAI_MOGYORODI}
 one can derive the following theorem.

\begin{Thm}\label{THM_MCLT_RAO2}
Let \ $\alpha\in\RR$.
\ Suppose that there exists a function \ $Q:[0,T)\to \RR \setminus \{ 0 \}$
 \ such that \ $\lim_{t\uparrow T}Q(t)=0$ \ and
 \begin{align}\label{SEGED_RAO}
   Q(t)^2\int_0^t\frac{b(s)^2(X_s^{(\alpha)})^2}{\sigma(s)^2}\,\dd s
      \stoch \eta^2
   \qquad \text{as}\quad t\uparrow T ,
 \end{align}
 where \ $\eta$ \ is a random variable defined on \ $(\Omega,\cF,\PP)$.
\ Then
 \begin{align*}
   \left(
      Q(t)\int_0^t\frac{b(s)X_s^{(\alpha)}}{\sigma(s)}\,\dd B_s\,,\,
          Q(t)^2\int_0^t\frac{b(s)^2(X_s^{(\alpha)})^2}{\sigma(s)^2}\,\dd s
   \right)
    \distr (\eta\xi,\eta^2)
 \qquad \text{as}\quad t\uparrow T ,
 \end{align*}
 where \ $\xi$ \ is a standard normally distributed random variable independent
 of \ $\eta$.
\ Moreover, if \ $\PP(\eta>0)=1$, \ then
 \begin{align*}
    \frac{1}{Q(t)}
         (\widehat\alpha_t^{(X^{(\alpha)})}-\alpha)
    \distr \frac{\xi}{\eta}
 \qquad \text{as}\quad t\uparrow T.
 \end{align*}
\end{Thm}

\noindent{\bf Proof.}
With the notation \ $M_t:=\int_0^t\frac{b(s)X_s^{(\alpha)}}{\sigma(s)}\,\dd B_s$, \ $t\in[0,T)$,
 \ we have \ $(M_t)_{t\in[0,T)}$ \ is a continuous square integrable martingale with respect to
 the filtration \ $(\cF_t)_{t\in[0,T)}$. \ By \eqref{SEGED_RAO} and Theorem \ref{THM_Zanten}, we have
 \[
   \left(
      Q(t)\int_0^t\frac{b(s)X_s^{(\alpha)}}{\sigma(s)}\,\dd B_s\,,\, \eta^2
   \right) \distr (\eta\xi,\eta^2)
   \qquad \text{as}\quad t\uparrow T .
 \]
By \eqref{SEGED_RAO} and Lemma \ref{LEMMA_KATAI_MOGYORODI}, we get
 \begin{align*}
     \left(
      Q(t)\int_0^t\frac{b(s)X_s^{(\alpha)}}{\sigma(s)}\,\dd B_s\,,\,
          Q(t)^2\int_0^t\frac{b(s)^2(X_s^{(\alpha)})^2}{\sigma(s)^2}\,\dd s
     \right)
   - \left(Q(t)\int_0^t\frac{b(s)X_s^{(\alpha)}}{\sigma(s)}\,\dd B_s\,,\,\eta^2 \right)
     \stoch 0,
 \end{align*}
 as \ $t\uparrow T$.
\ This implies the first part of the assertion using Slutsky's lemma (see, e.g.,
 van der Vaart \cite[Lemma 2.8]{Vaart}).
Using \eqref{SEGED20_uj} and the continuous mapping theorem (see, e.g.,
 van der Vaart \cite[Theorem 2.3]{Vaart}), we also have the second part of the assertion.
\proofend

Next we turn to the investigation of the asymptotic properties of the MLE of
 \ $\alpha$.

\begin{Thm} \label{sing_special}
Suppose that \ $\alpha \in \RR$ \ such that
 \begin{align}
  &\lim_{t \uparrow T} I_{X^{(\alpha)}}(t)=\infty,\label{Fisher5} \\
  &\lim_{t \uparrow T}
     \frac{b(t)}{\sigma(t)^2}
      \exp\left\{ 2 \alpha \int_0^t b(w) \, \dd w \right\}
      =C\in\RR\setminus\{0\}.\label{Fisher6}
 \end{align}
Then
 \begin{align*}%\label{sing_thm1}
   \sqrt{I_{X^{(\alpha)}}(t)} \, \Big( \halpha_t^{(X^{(\alpha)})} - \alpha \Big)
   \distr \frac{\sign(C)}{\sqrt{2}} \,
          \frac{\int_0^1 W_s \, \dd W_s}{\int_0^1 (W_s)^2 \, \dd s} \qquad
   \text{as} \quad t \uparrow T,
 \end{align*}
 where \ $\sign$ \ denotes the signum function and \ $(W_s)_{s\in [0,1]}$ \ is a standard Wiener process.
\end{Thm}

For the proof of Theorem \ref{sing_special} we need the following lemma.

\begin{Lem}\label{LEMMA_EKVIVALENS_SINGULAR}
Let \ $\alpha \in \RR$ \ be such that condition \eqref{Fisher6} is satisfied.
 Then \eqref{Fisher5} is equivalent to any of the following conditions:
 \begin{align}
   &\lim_{t \uparrow T}\int_0^t\vert b(v)\vert\,\dd v=\infty, \label{Fisher5_cor}\\
   &\lim_{t \uparrow T}\int_0^t\sigma(s)^2
      \,\exp\left\{ - 2 \alpha \int_0^s b(v)\,\dd v \right\} \dd s =\infty.
    \label{Fisher9}
 \end{align}
\end{Lem}

\noindent{\bf Proof.}
By \eqref{Fisher6}, there exist \ $c_1>0$, \ $c_2>0$ \ and \ $t_1\in(0,T)$ \ such
 that
 \begin{equation} \label{liminfsup}
   0<c_1 \vert b(t)\vert \exp\left\{ 2 \alpha \int_0^t b(w) \, \dd w \right\}
   \leq \sigma(t)^2
   \leq c_2 \vert b(t)\vert\exp\left\{ 2 \alpha \int_0^t b(w) \, \dd w \right\} ,
 \end{equation}
 for all \ $t \in [t_1,T)$.
\ First we show that \eqref{Fisher5} and \eqref{Fisher5_cor} are equivalent.
By \eqref{liminfsup}, we have for all \ $t\in[t_1,T)$, \
 \begin{align*}
  I_{X^{(\alpha)}}(t)-I_{X^{(\alpha)}}(t_1)
   &=\int_{t_1}^t\frac{b(s)^2}{\sigma(s)^2}
     \left(\int_0^s\sigma(u)^2
      \,\exp\left\{2\alpha\int_u^s b(v)\,\dd v \right\} \dd  u\right)\,\dd s\\[1mm]
     &\leq \int_{t_1}^t\frac{b(s)^2}{\sigma(s)^2}
      \exp\left\{2\alpha\int_0^s b(v)\,\dd v \right\}
      \left(\int_0^{t_1}\!\sigma(u)^2
      \,\exp\left\{-2\alpha\int_0^u b(v)\,\dd v \right\} \dd  u \right)\!\dd s\\[1mm]
     &\phantom{\leq\;}
       +\int_{t_1}^t\frac{b(s)^2}{\sigma(s)^2}
         \exp\left\{2\alpha\int_0^s b(v)\,\dd v \right\}
         \left(\int_{t_1}^sc_2 \vert b(u) \vert \, \dd  u\right)\,\dd s \\[1mm]
    &\leq a_1 \int_{t_1}^t \vert b(u) \vert \dd  u
          +a_2 \left( \int_{t_1}^t \vert b(u) \vert \, \dd  u \right)^2,
 \end{align*}
 where
 \begin{align*}
   a_1:=\frac{1}{c_1}
         \int_0^{t_1}\sigma(u)^2
           \,\exp\left\{-2\alpha\int_0^u b(v)\,\dd v \right\} \dd  u,
       \qquad \text{and}\qquad
  a_2:= \frac{c_2}{2c_1}.
 \end{align*}
Moreover, again by \eqref{liminfsup}, for all \ $t \in [t_1,T)$, \ we have
 \begin{align*}
  I_{X^{(\alpha)}}(t)-I_{X^{(\alpha)}}(t_1)
   \geq \frac{c_1}{2c_2}\left(\int_{t_1}^t \vert b(u)\vert \,\dd u\right)^2 .
 \end{align*}
This implies the equivalence of \eqref{Fisher5} and \eqref{Fisher5_cor}, since if
 \ $(x_n)_{n\in\NN}$ \ is a monotone increasing sequence of real numbers and
 \ $a_1>0$, \ $a_2>0$, \ then \ $a_1x_n+a_2x_n^2$ \ tends to
 \ $\infty$ \ if and only if \ $x_n\to\infty$.
\ Indeed, since \ $(x_n)_{n\in\NN}$ \ is monotone increasing,
 \ $\lim_{n\to\infty}x_n\in\RR$ \ exists or \ $x_n\uparrow\infty$.
\ In the first case  \ $a_1x_n+a_2x_n^2$ \  does not converge to \ $\infty.$

Now we show that \eqref{Fisher5_cor} and \eqref{Fisher9} are equivalent.
Using \eqref{liminfsup}, we have
 \begin{align*}
   \int_0^{t_1}\sigma(s)^2 \,& \exp\left\{ - 2 \alpha \int_0^s b(v) \, \dd v \right\} \dd s
            + c_1 \int_{t_1}^t \vert b(s) \vert \, \dd s\\
   &\leq \int_0^t\sigma(s)^2 \, \exp\left\{ - 2 \alpha \int_0^s b(v) \, \dd v \right\} \dd s\\
   &\leq \int_0^{t_1}\sigma(s)^2 \, \exp\left\{ - 2 \alpha \int_0^s b(v) \, \dd v \right\} \dd s
            + c_2 \int_{t_1}^t \vert b(s) \vert \, \dd s,\qquad t\in[t_1,T),
 \end{align*}
 which implies the corresponding part of the assertion.
\proofend

\noindent{\bf Proof of Theorem \ref{sing_special}.}
Note that condition \eqref{Fisher6} yields that there exists \ $t_0\in(0,T)$ \
 such that \ $b(t)\ne 0$ \ for all \ $t\in[t_0,T)$.
\ By Lemma \ref{LEMMA_EKVIVALENS_SINGULAR}, since \eqref{Fisher5} is assumed, we have conditions
 \eqref{Fisher5_cor} and \eqref{Fisher9} are also satisfied.
By \eqref{Fisher9}, for each
 \ $t \in (0,T)$, \ there exists a function \ $\tau_t : [0,\infty) \to [0,T)$
 \ such that
 \[
   \frac{1}{\sqrt{I_{X^{(\alpha)}}(t)}}
   \int_0^{\tau_t(u)}
    \sigma(s)^2 \, \exp\left\{ - 2 \alpha \int_0^s b(w) \, \dd w \right\} \dd s
   = u, \qquad \text{$u \in [0,\infty)$.}
 \]
Clearly, \ $\tau_t$ \ is strictly increasing (hence invertible), and again by \eqref{Fisher9},
 \ $\lim\limits_{u \to \infty} \tau_t(u) = T$ \ for all \ $t\in(0,T)$, \ and
 \[
   \tau^{-1}_t(v)
   = \frac{1}{\sqrt{I_{X^{(\alpha)}}(t)}}
     \int_0^v
      \sigma(s)^2 \,
      \exp\left\{ - 2 \alpha \int_0^s b(w) \, \dd w \right\} \dd s, \qquad
   \text{$v \in [0,T)$.}
 \]
Then \ $\lim\limits_{v \uparrow T} \tau_t^{-1}(v) = \infty$ \ for all \ $t\in(0,T)$, \  and
 \begin{align*} %\label{SEGED47}
     (\tau^{-1}_t)'(v)
   = \frac{1}{\sqrt{I_{X^{(\alpha)}}(t)}} \,
     \sigma(v)^2 \,
     \exp\left\{ - 2 \alpha \int_0^v b(w) \, \dd w \right\} , \qquad
   \text{$v \in (0,T)$.}
 \end{align*}
By the theorem on differentiation of inverse function, \ $\tau_t$ \ is also continuously differentiable
 and
 \[
  (\tau_t)'(u)
    = \sqrt{I_{X^{(\alpha)}}(t)}\sigma(\tau_t(u))^{-2}
      \exp\left\{ 2 \alpha \int_0^{\tau_t(u)} b(v) \, \dd v \right\},
      \qquad u\in(0,\infty).
 \]
The process
 \begin{align} \label{M}
  \begin{split}
    M_t^{(X^{(\alpha)})}
    &:= X_t^{(\alpha)} \exp\left\{ - \alpha \int_0^t b(u) \, \dd u \right\}\\
    & = \int_0^t
       \sigma(s)
       \exp\left\{ - \alpha \int_0^s b(u) \, \dd u \right\} \, \dd B_s ,
    \qquad t \in [0,T) ,
   \end{split}
 \end{align}
 is a continuous square-integrable martingale with respect to the filtration
 induced by \ $B$.
\ With this notation we have for all \ $t\in(0,T)$,
\begin{align*}%\label{SEGED_SINGULAR4}
  \frac{1}{I_{X^{(\alpha)}}(t)}
  \int_0^t
    &\frac{b(v)^2 \, ( X_v^{(\alpha)} )^2}{\sigma(v)^2} \, \dd v \\
    &= \frac{1}{I_{X^{(\alpha)}}(t)}
     \int_0^t
      \frac{b(v)^2}{\sigma(v)^2} \,
      \exp\left\{ 2 \alpha \int_0^v b(w) \, \dd w \right\}
      \big( M_v^{(X^{(\alpha)})} \big)^2 \, \dd v \\[1mm]
   &= \frac{1}{\sqrt{I_{X^{(\alpha)}}(t)}}
     \int_0^{\tau_t^{-1}(t)}
      \frac{b(\tau_t(u))^2}{\sigma(\tau_t(u))^4} \,
      \exp\left\{ 4 \alpha \int_0^{\tau_t(u)} b(w) \, \dd w \right\}
      \big( M_{\tau_t(u)}^{(X^{(\alpha)})} \big)^2 \, \dd u.
\end{align*}
Then for all \ $t\in(0,T)$,
 \begin{align}\label{SEGED_SINGULAR4}
  \frac{1}{I_{X^{(\alpha)}}(t)}
  \int_0^t
    \frac{b(v)^2 \, ( X_v^{(\alpha)} )^2}{\sigma(v)^2} \, \dd v
   = \int_0^{\tau_t^{-1}(t)}
      c(\tau_t(u))^2 \,
      \big( \tM_u^{(X^{(\alpha)},\,t)} \big)^2 \, \dd u ,
\end{align}
 where
 \[
   c(s):=\frac{b(s)}{\sigma(s)^2}
           \exp\left\{ 2 \alpha \int_0^s b(w) \, \dd w \right\},
   \qquad s\in[0,T),
 \]
 and
 \[
   \tM_u^{(X^{(\alpha)},\,t)}
   := \frac{1}{\sqrt[4]{I_{X^{(\alpha)}}(t)}} M_{\tau_t(u)}^{(X^{(\alpha)})} ,
   \qquad
   u \in [0,\infty).
 \]
By \eqref{Fisher6}, we have \ $\lim_{s\uparrow T}c(s)=C$,
 \ and for all \ $t\in(0,T)$, \ the process \ $(\tM_u^{(X^{(\alpha)},\,t)})_{u\in[0,\infty)}$
 \ is a continuous Gauss martingale with respect to the filtration \ $(\widetilde\cF_u^t)_{u\geq 0}$,
 \ where
 \[
  \widetilde\cF_u^t:=\sigma\big(B_v, 0\leq v\leq \tau_t(u)\big), \qquad u\geq 0.
  \]
Moreover, for all \ $t\in(0,T)$, \ the process \ $(\tM_u^{(X^{(\alpha)},\,t)})_{u\in[0,\infty)}$ \ has
 quadratic variation
 \[
   \langle \tM^{(X^{(\alpha)},\,t)} \rangle_u
   = \frac{1}{\sqrt{I_{X^{(\alpha)}}(t)}}
     \int_0^{\tau_t(u)}
      \sigma(s)^2 \,
      \exp\left\{ - 2 \alpha \int_0^s b(w) \, \dd w \right\} \dd s
   = u , \qquad  u \in [0,\infty).
 \]
Then Theorem 3.3.16 in Karatzas and Shreve \cite{KarShr} yields that
 \ $(\tM_u^{(X^{(\alpha)},\,t)})_{u \in [0,\infty)}$ \ is a standard Wiener process
 with respect to the filtration  \ $(\widetilde\cF_u^t)_{u\geq 0}$.
\ In a similar way we get
 \begin{align}\label{SEGED_SINGULAR5}
   \begin{split}
     \frac{1}{\sqrt{I_{X^{(\alpha)}}(t)}}
     &\int_0^t \frac{b(v) \, X_v^{(\alpha) }}{\sigma(v)} \, \dd B_v \\
     &= \frac{1}{\sqrt{I_{X^{(\alpha)}}(t)}}
        \int_0^t
          \frac{b(v)}{\sigma(v)^2} \,
          \exp\left\{ 2 \alpha \int_0^v b(w) \, \dd w \right\}
          M_v^{(X^{(\alpha)})} \, \dd M_v^{(X^{(\alpha)})} \\[1mm]
     &= \frac{1}{\sqrt{I_{X^{(\alpha)}}(t)}}
         \int_0^{\tau_t^{-1}(t)}
         \frac{b(\tau_t(u))}{\sigma(\tau_t(u))^2} \,
         \exp\left\{ 2 \alpha \int_0^{\tau_t(u)} b(w) \, \dd w \right\}
         M_{\tau_t(u)}^{(X^{(\alpha)})} \, \dd M_{\tau_t(u)}^{(X^{(\alpha)})} \\[1mm]
     &= \int_0^{\tau_t^{-1}(t)}
        c(\tau_t(u)) \, \tM_u^{(X^{(\alpha)},\,t)} \, \dd \tM_u^{(X^{(\alpha)},\,t)} ,
        \qquad t\in(0,T),
   \end{split}
 \end{align}
 where the last but one equality follows by the construction of a stochastic integral
 with respect to \ $M^{(X^{(\alpha)})}$, \  see, e.g.,
 Jacod and Shiryaev \cite[Proposition 4.44, Chapter I]{JacShi}.
By assumption \eqref{Fisher5} and the fact that \ $b(t)\ne0$ \ for all \ $t\in[t_0,T)$,
 \ we can use L'Hospital's rule and we obtain
 \begin{align}\label{SEGED_SINGULAR6}
   \begin{split}
    \lim_{t \uparrow T} \big( \tau_t^{-1}(t) \big)^2
     &=\lim_{t \uparrow T}
       \frac{\left(\int_0^t
                  \sigma(s)^2 \,
                  \exp\left\{ - 2 \alpha \int_0^s b(w) \, \dd w \right\}
                  \dd s \right)^2}
          {I_{X^{(\alpha)}}(t)} \\[1mm]
     &=\lim_{t \uparrow T}
        \frac{2 \sigma(t)^2 \,
           \exp\left\{ - 2 \alpha \int_0^t b(w) \, \dd w \right\}
           \int_0^t
            \sigma(s)^2 \,
            \exp\left\{ - 2 \alpha \int_0^s b(w) \, \dd w \right\} \dd s}
          {\frac{b(t)^2}{\sigma(t)^2}
           \int_0^t
            \sigma(s)^2 \,
            \exp\left\{ 2 \alpha \int_s^t b(w) \, \dd w \right\} \dd s} \\[1mm]
    &=\lim_{t \uparrow T}
       \frac{2 \, \sigma(t)^4}{b(t)^2}
      \exp\left\{ - 4 \alpha \int_0^t b(w) \, \dd w \right\}
     = \lim_{t \uparrow T}\frac{2}{c(t)^2}
     = \frac{2}{C^2},
  \end{split}
 \end{align}
 where the last equality follows by \eqref{Fisher6}.
Hence, using that \ $\tau_t^{-1}(t)\in[0,\infty)$, \ we also have
 \ $\lim_{t \uparrow T}\tau_t^{-1}(t)=\frac{\sqrt{2}}{\vert C\vert}$.

Now we prove that
 \begin{align} \label{joint}
  \begin{split}
   &\left( \frac{1}{\sqrt{I_{X^{(\alpha)}}(t)}}
                 \int_0^t \frac{b(s) \, X_s^{(\alpha)}}{\sigma(s)} \, \dd B_s
               , \,
          \frac{1}{I_{X^{(\alpha)}}(t)}
              \int_0^t \frac{b(s)^2 \, \left(X_s^{(\alpha)}\right)^2}
                             {\sigma(s)^2} \, \dd s
                \right)\\[1mm]
   &\qquad\qquad
    \distr \left( C \int_0^{\sqrt{2} / |C|} W_s \, \dd W_s, \,
                 C^2 \int_0^{\sqrt{2} / |C|} (W_s)^2 \, \dd s \right)
         \qquad \text{as \ $t \uparrow T$.}
  \end{split}
 \end{align}
Using \eqref{SEGED_SINGULAR4}, \eqref{SEGED_SINGULAR5} and that
 \ $(\tM_u^{(X^{(\alpha)},\,t)})_{u \in [0,\infty)}$ \ is a standard Wiener process
 for all \ $t\in(0,T)$, \ we conclude that
 \[
  \left( \frac{1}{\sqrt{I_{X^{(\alpha)}}(t)}}
                 \int_0^t \frac{b(s) \, X_s^{(\alpha)}}{\sigma(s)} \, \dd B_s
               , \,
          \frac{1}{I_{X^{(\alpha)}}(t)}
              \int_0^t \frac{b(s)^2 \, \left(X_s^{(\alpha)}\right)^2}
                             {\sigma(s)^2} \, \dd s
                \right)
 \]
 has the same distribution as
 \[
  \left( \int_0^{\tau_t^{-1}(t)} c(\tau_t(u)) \, W_u \, \dd W_u, \,
          \int_0^{\tau_t^{-1}(t)} c(\tau_t(u))^2 \, (W_u)^2 \, \dd u
    \right),
 \]
 for all \ $t\in(0,T)$ \ with some fixed standard Wiener process
 \ $(W_u)_{u\geq 0}$.
\ Hence to prove \eqref{joint}, using Slutsky's lemma, it is enough to check
 that
 \begin{align*}
  &\left( \int_0^{\tau_t^{-1}(t)} c(\tau_t(u)) \, W_u \, \dd W_u, \,
          \int_0^{\tau_t^{-1}(t)} c(\tau_t(u))^2 \, (W_u)^2 \, \dd u
    \right) \\
  & \qquad
    -\left( C \int_0^{\sqrt{2} / |C|} W_u \, \dd W_u, \,
            C^2 \int_0^{\sqrt{2} / |C|} (W_u)^2 \, \dd u \right)
  \stoch 0
  \qquad \text{as \ $t \uparrow T$.}
 \end{align*}
For this it is enough to prove that the following convergences hold:
 \begin{align}\label{c1}
    &\int_0^{\sqrt{2} / |C|} \big[ c(\tau_t(u)) - C \big] W_u \, \dd W_u
        \stackrel{L_2}{\longrightarrow} 0
          \qquad \text{as \ $t \uparrow T,$}\\[1mm]\label{c2}
    &\int_0^{\tau_t^{-1}(t)} c(\tau_t(u)) \, W_u \, \dd W_u
        -\int_0^{\sqrt{2} / |C|} c(\tau_t(u)) \, W_u \, \dd W_u
          \stackrel{L_2}{\longrightarrow} 0
       \qquad \text{as \ $t \uparrow T,$}\\[1mm]\label{c3}
    &\PP\left(\lim_{t\uparrow T}
               \int_0^{\sqrt{2} / |C|} \big[ c(\tau_t(u))^2 - C^2 \big] \, (W_u)^2 \, \dd u
               =0\right)=1,\\[1mm]\label{c4}
    &\PP\left(\lim_{t\uparrow T}
               \left(\int_0^{\tau_t^{-1}(t)} c(\tau_t(u))^2 \, (W_u)^2 \, \dd u
               -\int_0^{\sqrt{2} / |C|} c(\tau_t(u))^2 \, (W_u)^2 \, \dd u\right)
               =0\right)=1.
 \end{align}
Using that \ $\lim\limits_{v \to \infty} \tau_t(v) = T$ \ for all
 \ $t\in(0,T)$ \ and \ $\lim\limits_{t \uparrow T} \tau_t(v) = T$ \
 for all \ $v\in(0,\infty),$ \ first we prove \eqref{c1}.
An easy calculation shows that for all \ $t\in(0,T)$,
 \begin{multline*}
  \EE\left(\int_0^{\sqrt{2} / |C|}
            \big[ c(\tau_t(u)) - C \big] W_u \, \dd W_u \right)^2
  = \EE\int_0^{\sqrt{2} / |C|} \big[ c(\tau_t(u)) - C \big]^2 (W_u)^2 \, \dd u \\
  = \int_0^{\sqrt{2} / |C|} \big[ c(\tau_t(u)) - C \big]^2  u \, \dd u
  \leq \frac{\sqrt{2}}{|C|}
       \int_0^{\sqrt{2} / |C|} \big[ c(\tau_t(u)) - C \big]^2 \dd u .
 \end{multline*}
The only non-trivial step is to verify that the first equality holds.
By Karatzas and Shreve \cite[Proposition 3.2.10]{KarShr}, for this equality
 it is enough to check that
 \[
    \EE\int_0^{\sqrt{2} / |C|} \big[ c(\tau_t(u)) - C \big]^2 (W_u)^2 \, \dd u
     = \int_0^{\sqrt{2} / |C|} \big[ c(\tau_t(u)) - C \big]^2  u \, \dd u
       < \infty, \qquad t\in(0,T) ,
 \]
 which holds, since the integrand
 \ $u\mapsto \big[ c(\tau_t(u)) - C \big]^2 u$ \ is continuous on \ $[0,\sqrt{2} / |C|]$ \
 and hence bounded.
Finally, we prove that
 \begin{align}\label{SEGED_SINGULAR9}
   \lim_{t\uparrow T}
     \int_{0}^{\sqrt{2} / |C|} \big[ c(\tau_t(u)) - C \big]^2 \dd u
      =0.
 \end{align}
Since \ $\lim_{s\uparrow T}c(s)=C\in\RR\setminus\{0\}$, \
 for all \ $\delta > 0$ \ there exists \ $\vare > 0$ \ such that
 \ $|c(s) - C| < \delta$ \ for all \ $s \in (T - \vare, T)$.
\ For all \ $\vare>0$ \ and for all \ $u_0 \in (0,\infty)$ \ there exists
 \ $t_1 \in (0,T)$ \ such that \ $\tau_t(u_0) \in (T - \vare, T)$ \ for all
 \ $t \in (t_1,T)$, \ and hence \ $\tau_t(u) \in (T - \vare, T)$ \ for all
 \ $t \in (t_1,T)$ \ and \ $u \geq u_0$, \ since \ $\tau_t$ \ is increasing.
Consequently, for all \ $\delta > 0$ \ and all \ $u_0 \in (0,\infty)$ \ there
 exists \ $t_1 \in (0,T)$ \ such that \ $|c(\tau_t(u)) - C| < \delta$ \ for all
 \ $t \in (t_1,T)$ \ and all \ $u \geq u_0$. \
Thus for all \ $\delta>0$ \ and all \ $u_0\in(0,\sqrt{2} \big/\, |C|)$, \ there exists
 \ $t_1\in(0,T)$ \ such that
 \[
   \int_{u_0}^{\sqrt{2} / |C|} \big[ c(\tau_t(u)) - C \big]^2 \dd u
   \leq \frac{\sqrt{2}}{|C|} \, \delta^2,
   \qquad\text{$t\in(t_1,T)$.}
 \]
Then for all \ $\delta>0$ \ and all \ $u_0\in(0,\sqrt{2} \big/\, |C|)$, \ there exists
 \ $t_1\in(0,T)$ \ such that
 \begin{align*}
   \int_{0}^{\sqrt{2} / |C|} \big[ c(\tau_t(u)) - C \big]^2 \dd u
       &=\int_0^{u_0} \big[ c(\tau_t(u)) - C \big]^2 \dd u
        +\int_{u_0}^{\sqrt{2} / |C|} \big[ c(\tau_t(u)) - C \big]^2 \dd u \\
       &\leq \sup_{(u,t)\in[0,u_0]\times[t_1,T)}\big[ c(\tau_t(u)) - C \big]^2u_0
             +\frac{\sqrt{2}}{|C|} \, \delta^2,
        \qquad \text{$t\in(t_1,T)$.}
 \end{align*}
Since \ $\lim_{s\uparrow T}c(s)=C\in\RR\setminus\{0\}$ \ implies that there exists
 \ $K_1\in(0,\infty)$ \ such that
 \[
    \sup_{(u,t)\in[0,u_0]\times[t_1,T)}\big[ c(\tau_t(u)) - C \big]^2
         \leq \sup_{s\in[0,T)}(c(s) - C)^2\leq K_1,
 \]
 we have for all \ $\delta>0$ \ and all \ $u_0\in(0,\sqrt{2} \big/\, |C|)$, \ there exists
 \ $t_1\in(0,T)$ \ such that
 \[
    \int_{0}^{\sqrt{2} / |C|} \big[ c(\tau_t(u)) - C \big]^2 \dd u
       \leq K_1u_0+\frac{\sqrt{2}}{|C|}\delta^2,
       \qquad \text{$t\in(t_1,T)$,}
 \]
 which yields \eqref{SEGED_SINGULAR9} and then we obtain \eqref{c1}.

Now we check \eqref{c2}.
Similarly as above, we have for all \ $t\in(0,T)$,
 \begin{align*}
   &\EE\left(
          \int_0^{\tau_t^{-1}(t)} c(\tau_t(u)) \, W_u \, \dd W_u
          -\int_0^{\sqrt{2} / |C|} c(\tau_t(u)) \, W_u \, \dd W_u
       \right)^2
    =\int_{\frac{\sqrt{2}}{|C|}\wedge \tau_t^{-1}(t)}^{\frac{\sqrt{2}}{|C|}\vee \tau_t^{-1}(t)}
           c(\tau_t(u))^2u\,\dd u\\[1mm]
   &\phantom{\qquad\qquad}
    \leq K_2\left(\frac{\sqrt{2}}{|C|}\vee \tau_t^{-1}(t)\right)
          \left\vert \frac{\sqrt{2}}{|C|} - \tau_t^{-1}(t) \right\vert ,
 \end{align*}
 where the last step follows from \ $K_2:=\sup_{s\in[0,T)}c(s)^2<\infty$, \
 since \ $\lim_{s\uparrow T}c(s)=C\in\RR\setminus\{0\}$.
\ By \eqref{SEGED_SINGULAR6}, we have \ $\lim_{t\uparrow T}\tau_t^{-1}(t)=\sqrt{2}/|C|$ \
 and hence
 \[
    \lim_{t\uparrow T}
       \left\vert \frac{\sqrt{2}}{|C|} - \tau_t^{-1}(t) \right\vert
       =0,
 \]
 which implies \eqref{c2}.

Now we check \eqref{c3}.
Using that, by Cauchy-Schwartz's inequality,
 \begin{align*}
   &\left(\int_0^{\sqrt{2} / |C|}
           \big[ c(\tau_t(u))^2 - C^2 \big] \, (W_u)^2 \, \dd u \right)^2\\
     &\phantom{\qquad\qquad}
      \leq \left( \int_0^{\sqrt{2} / |C|}
                \big[ c(\tau_t(u))^2 - C^2 \big]^2 \, \dd u \right)
        \left( \int_0^{\sqrt{2} / |C|} (W_u)^4 \, \dd u \right),
         \qquad \text{$t\in(0,T),$}
 \end{align*}
 we have it is enough to check that
 \[
   \lim_{t\uparrow T}
     \int_{0}^{\sqrt{2} / |C|} \big[ c(\tau_t(u))^2 - C^2 \big]^2 \dd u
      =0,
 \]
 since
 \[
   \PP\left(\int_{0}^{\sqrt{2}/|C|}(W_u)^4\,\dd u<\infty\right)=1.
 \]
Using that \ $\lim_{s\uparrow T}c(s)=C\in\RR\setminus\{0\}$ \ and that \ $c$ \ is
 continuous, there exists \ $K_3\in(0,\infty)$ \ such that
 \ $\sup_{s\in[0,T)}\vert c(s)\vert\leq K_3$.
\ Hence
 \begin{align*}
   \int_{0}^{\sqrt{2} / |C|} \big[ c(\tau_t(u))^2 - C^2 \big]^2 \dd u
     &=\int_{0}^{\sqrt{2} / |C|} \big[ c(\tau_t(u)) + C \big]^2\big[ c(\tau_t(u)) - C \big]^2 \dd u\\
     &\leq (K_3+\vert C\vert)^2
             \int_{0}^{\sqrt{2} / |C|}\big[ c(\tau_t(u)) - C \big]^2 \dd u
     \to 0 \qquad \text{as \ $t\uparrow T$,}
 \end{align*}
 where the last step follows by \eqref{SEGED_SINGULAR9}.

Using the very same arguments as above, one can check \eqref{c4}.

By \eqref{joint} and the continuous mapping theorem, we have
 \begin{align*}
  \sqrt{I_{X^{(\alpha)}}(t)} \, ( \widehat\alpha_t^{(X^{(\alpha)})} - \alpha )
  = \frac{\frac{1}{\sqrt{I_{X^{(\alpha)}}(t)}}
           \int_0^t \frac{b(s) \, X_s^{(\alpha)}}{\sigma(s)} \, \dd B_s}
          {\frac{1}{I_{X^{(\alpha)}}(t)}
          \int_0^t \frac{b(s)^2 \, ( X_s^{(\alpha)} )^2}{\sigma(s)^2} \, \dd s}
  \distr \frac{C \int_0^{\sqrt{2} / |C|} W_s \, \dd W_s}
               {C^2 \int_0^{\sqrt{2} / |C|} (W_s)^2 \, \dd s}
   \qquad \text{as} \quad  t \uparrow T.
 \end{align*}
Using that for all \ $\lambda>0$, \ the process
 \ $\left(\lambda^{-1/2}W_{\lambda t}\right)_{t\geq 0}$  \ is a standard Wiener process, by the
 substitution \ $s=\frac{\sqrt 2}{\vert C\vert}u$, \ $u\in\RR$, \ we get the random variable
 \begin{align*}
  \frac{C \int_0^{\sqrt{2} / |C|} W_s \, \dd W_s}
               {C^2 \int_0^{\sqrt{2} / |C|} (W_s)^2 \, \dd s}
  =
  \frac{1}{C}
    \frac{\int_0^1 W_{\frac{\sqrt 2}{\vert C\vert}u} \, \dd W_{\frac{\sqrt 2}{\vert C\vert}u}}
               {\int_0^1 (W_{\frac{\sqrt 2}{\vert C\vert}u})^2 \frac{\sqrt 2}{\vert C\vert}\, \dd u}
 =\frac{1}{C}
     \frac{\frac{\sqrt 2}{\vert C\vert}
        \int_0^1 \sqrt{\frac{\vert C\vert}{\sqrt 2}} W_{\frac{\sqrt 2}{\vert C\vert}u}
                  \, \dd \sqrt{\frac{\vert C\vert}{\sqrt 2}} W_{\frac{\sqrt 2}{\vert C\vert}u}}
     {\frac{\sqrt 2}{\vert C\vert} \frac{\sqrt 2}{\vert C\vert}
        \int_0^1\Big(\sqrt{\frac{\vert C\vert}{\sqrt 2}}W_{\frac{\sqrt 2}{\vert C\vert}u}\Big)^2
        \,\dd u}
 \end{align*}
 has the same distribution as
 \[
    \frac{\sign(C)}{\sqrt{2}} \,
           \frac{\int_0^1 W_s \, \dd W_s}{\int_0^1 (W_s)^2 \, \dd s}.
 \]
\proofend

For historical fidelity, we remark that the corresponding part of Example 8.1 in Luschgy \cite{Lus2}
 is a special case of our Theorem \ref{sing_special}, and in our proof we used some ideas
 of Luschgy's example.
Note also that, by Lemma \ref{LEMMA_EKVIVALENS_SINGULAR}, condition \eqref{Fisher5}
 in Theorem \ref{sing_special} can be replaced by \eqref{Fisher5_cor} or \eqref{Fisher9}.

In the next remark we give an example for functions \ $b$ \ and \ $\sigma$ \ for which
 conditions \eqref{Fisher5} and \eqref{Fisher6} are satisfied.

\begin{Rem}
First let \ $\alpha\ne0$.
\ Let \ $\sigma:[0,T)\to(0,\infty)$ \ be some continuously differentiable function such that
 \ $\int_0^T\sigma(s)^2\,\dd s:=\lim_{u\uparrow T}\int_0^u\sigma(s)^2\,\dd s<\infty$, \ and
 let
 \[
    b(t):=-\frac{1}{2\alpha}\frac{\sigma(t)^2}{\int_t^T\sigma(s)^2\,\dd s},
         \qquad t\in[0,T).
 \]
Then for all \ $0\leq s<t<T$,
 \begin{align*}
    \int_s^t \vert b(u)\vert\,\dd u
        =-\frac{1}{2\vert\alpha\vert}
           \ln\left(\frac{\int_t^T\sigma(v)^2\,\dd v}{\int_s^T\sigma(v)^2\,\dd v}\right),
 \end{align*}
 which implies that \ $\lim_{t\uparrow T}\int_0^t \vert b(u)\vert\,\dd u=\infty$.
\ Moreover,
 \begin{align*}
   \frac{b(t)}{\sigma(t)^2}\exp\left\{2\alpha\int_0^tb(u)\,\dd u\right\}
        =-\frac{1}{2 \alpha}\frac{1}{\int_0^T\sigma(s)^2\,\dd s}\in\RR\setminus\{0\},
 \end{align*}
 which implies \eqref{Fisher6}.
By Lemma \ref{LEMMA_EKVIVALENS_SINGULAR}, since
 \ $\lim_{t\uparrow T}\int_0^t \vert b(u)\vert\,\dd u=\infty$,
 \ we have \eqref{Fisher5} is also satisfied.

\noindent Let us suppose now that \ $\alpha=0$.
\ Let \ $\sigma:[0,T)\to(0,\infty)$ \ be some continuously differentiable function such that
 \ $\int_0^T\sigma(s)^2\,\dd s=\infty$, \ and \ $b(t):=\sigma(t)^2$, \ $t\in[0,T)$.
\ Then we have \ $\lim_{t\uparrow T}\int_0^t \vert b(u)\vert\,\dd u=\infty$, \ and,
 with \ $\alpha=0$,
 \[
   \frac{b(t)}{\sigma(t)^2}\exp\left\{2\alpha\int_0^tb(u)\,\dd u\right\}
      =1\in\RR\setminus\{0\},\qquad t\in[0,T),
 \]
 which implies \eqref{Fisher6}.
By Lemma \ref{LEMMA_EKVIVALENS_SINGULAR}, since
 \ $\lim_{t\uparrow T}\int_0^t \vert b(u)\vert\,\dd u=\infty$, \
 we have \eqref{Fisher5} is also satisfied.
\end{Rem}

Next we deal with the case of Cauchy limit distribution.

\begin{Thm}\label{PRO_CAUCHY}
Suppose that \ $\alpha \in \RR$ \ such that
 \begin{gather}
  \lim_{t \uparrow T} I_{X^{(\alpha)}}(t) = \infty , \label{Fisher1} \\
  \lim_{t\uparrow T}\int_0^t
     \sigma(s)^2 \exp\left\{ - 2 \alpha \int_0^s b(v) \, \dd v \right\} \, \dd s
      < \infty . \label{Fisher2}
 \end{gather}
Then
 \[
   \sqrt{I_{X^{(\alpha)}}(t)} \, ( \halpha_t^{(X^{(\alpha)})} - \alpha )
   \distr \zeta \qquad
   \text{as \ $t \uparrow T$,}
 \]
 where \ $\zeta$ \ is a random variable with standard Cauchy distribution
 admitting a density function \ $\frac{1}{\pi(1+x^2)}$, \ $x\in\RR$.
\end{Thm}

\noindent{\bf Proof.}
The process \ $(M_t^{(X^{(\alpha)})})_{t \in [0,T)}$ \ introduced in \eqref{M} is a
 continuous square-integrable martingale with respect to the filtration
 induced by \ $B$ \ and with quadratic variation
 \[
   \langle M^{(X^{(\alpha)})} \rangle_t
   = \int_0^t
      \sigma(s)^2
      \exp\left\{ - 2 \alpha \int_0^s b(u) \, \dd u \right\} \, \dd s,
   \qquad t \in [0,T) .
 \]
By \eqref{Fisher2}, we have
 \ $\lim_{t \uparrow T} \, \langle M^{(X^{(\alpha)})} \rangle_t < \infty$.
\ Hence Proposition 1.26 in Chapter IV and Proposition 1.8 in Chapter V in Revuz
 and Yor \cite{RevYor} imply that the limit
 \ $M^{(X^{(\alpha)})}_T := \lim_{ t\uparrow T} M^{(X^{(\alpha)})}_t$ \ exists almost surely.
Since \ $M^{(X^{(\alpha)})}_t$ \ is normally distributed with mean \ $0$ \ and with
 variance \ $\langle M^{(X^{(\alpha)})} \rangle_t$ \ for all \ $t \in [0,T)$, \ the
 random variable \ $M^{(X^{(\alpha)})}_T$ \ is also normally distributed with mean
 \ $0$ \ and with variance \ $\lim_{t\uparrow T} \langle M^{(X^{(\alpha)})} \rangle_t$.
\ Indeed, normally distributed random variables can converge in distribution
 only to a normally distributed random variable, by continuity theorem, see,
 e.g., page 304 in Shiryaev \cite{Shi}.
Hence for the random variable \ $M^{(X^{(\alpha)})}_T$ \ we obtain
 \[
   \PP\left( \lim_{t \uparrow T} M^{(X^{(\alpha)})}_t = M^{(X^{(\alpha)})}_T \right) = 1
   \qquad \text{and} \qquad
   M^{(X^{(\alpha)})}_T
   \distre \cN \left(0, \lim_{t\uparrow T} \langle M^{(X^{(\alpha)})} \rangle_t \right).
 \]
By \eqref{Fisher1} and the fact that \ $b(t)\ne0$ \ for all \ $t\in[t_0,T)$,
 \ we can use L'Hospital's rule and we obtain
 \begin{align*}
  \lim_{t \uparrow T}
   &\frac{\int_0^t
          \frac{b(s)^2 \, ( X_s^{(\alpha)} )^2}{\sigma(s)^2} \, \dd s}
        {I_{X^{(\alpha)}}(t)}
   = \lim_{t \uparrow T}
   \frac{\int_0^t
          \frac{b(s)^2 \, ( X_s^{(\alpha)} )^2}{\sigma(s)^2} \, \dd s}
        {\int_0^t
         \frac{b(s)^2 \, \EE ( X_s^{(\alpha)} )^2}{\sigma(s)^2} \, \dd s}
  = \lim_{t \uparrow T}
     \frac{( X_t^{(\alpha)} )^2}{\EE ( X_t^{(\alpha)} )^2} \\[2mm]
  &= \lim_{t \uparrow T}
     \frac{( M_t^{(X^{(\alpha)})} )^2
           \exp\left\{ 2 \alpha \int_0^t b(v) \, \dd v \right\}}
          {\int_0^t
            \sigma(u)^2
            \exp\left\{ 2 \alpha \int_u^t b(v) \, \dd v \right\} \, \dd u}
  = \lim_{t \uparrow T}
     \frac{( M_t^{(X^{(\alpha)})} )^2}
          {\int_0^t
            \sigma(u)^2
            \exp\left\{ - 2 \alpha \int_0^u b(v) \, \dd v \right\} \, \dd u}
    \\[2mm]
  &= \frac{( M_T^{(X^{(\alpha)})} )^2}
         {\lim_{t \uparrow T}\langle M^{(X^{(\alpha)})} \rangle_t}
  =: \xi^2
  \qquad \text{$\PP$-almost surely,}
 \end{align*}
 where \ $\xi \distre \cN(0,1)$.
\ Using Theorem \ref{THM_MCLT_RAO2} with \ $Q(t):=1/\sqrt{I_{X^{(\alpha)}}(t)}$,
 \ $t\in(0,T)$, \ we have
 \[
   \sqrt{I_{X^{(\alpha)}}(t)} \, ( \halpha_t^{(X^{(\alpha)})} - \alpha )
      \distr \frac{\eta}{|\xi|}.
 \]
This yields the assertion, since one can easily check that \ $\frac{\eta}{|\xi|}$ \ has
 standard Cauchy distribution.
\proofend

\begin{Rem}\label{REM_LUSCHGY1}
We note that if condition \eqref{Fisher2} is satisfied then
 \ $\lim_{t\uparrow T}\int_0^t\vert b(s)\vert\,\dd s=\infty$ \ yields condition \eqref{Fisher1}.
Indeed, for all \ $t\in(0,T)$,
 \begin{align*}
   I_{X^{(\alpha)}}(t)
       =\int_0^t\frac{b(s)^2}{\sigma(s)^2}
           \exp\left\{2\alpha\int_0^s b(u)\,\dd u\right\}
         \int_0^s\sigma(u)^2\exp\left\{-2\alpha\int_0^u b(v)\,\dd v\right\}\,\dd u\,\dd s,
 \end{align*}
 and, by \eqref{Fisher2}, for all \ $\varepsilon>0$ \ there exists a \ $t_\varepsilon\in[0,T)$
 \ such that for all \ $t\in[t_\varepsilon,T)$,
 \[
   \left\vert\int_0^t\sigma(s)^2\exp\left\{-2\alpha\int_0^s b(v)\,\dd v\right\}\,\dd s
              - \int_0^T\sigma(s)^2\exp\left\{-2\alpha\int_0^s b(v)\,\dd v\right\}\,\dd s
    \right\vert<\varepsilon.
 \]
Hence for all
 \[
   0<\varepsilon<\int_0^T\sigma(s)^2\exp\left\{-2\alpha\int_0^s b(v)\,\dd v\right\}\,\dd s,
 \]
 and for all \ $t\in[t_\varepsilon,T)$, \ we have
 \begin{align*}
   I_{X^{(\alpha)}}(t)
     \geq&\int_0^{t_\varepsilon}\frac{b(s)^2}{\sigma(s)^2}
          \int_0^s\sigma(u)^2\exp\left\{2\alpha\int_u^s b(v)\,\dd v\right\}\,\dd u\,\dd s\\[1mm]
         &+\int_{t_\varepsilon}^t\frac{b(s)^2}{\sigma(s)^2}
          \exp\left\{2\alpha\int_0^s b(u)\,\dd u\right\}\,\dd s\,
          \left(\int_0^T\sigma(s)^2\exp\left\{-2\alpha\int_0^s b(v)\,\dd v\right\}\,\dd s-\varepsilon\right).
 \end{align*}
This yields that it is enough to check that
 \begin{align}\label{SEGED_LUSCHGY1}
    \lim_{t\uparrow T}
       \int_0^t\frac{b(s)^2}{\sigma(s)^2}
         \exp\left\{2\alpha\int_0^s b(u)\,\dd u\right\}\,\dd s
       =\infty.
 \end{align}
By Cauchy-Schwartz's inequality, we get for all \ $t\in[0,T)$,
 \begin{align*}
    \left(\int_0^t\vert b(s)\vert\,\dd s\right)^2
     &=\left(\int_0^t\frac{\vert b(s)\vert}{\sigma(s)}
                \exp\left\{\alpha\int_0^s b(u)\,\dd u\right\}
                \sigma(s)
                \exp\left\{-\alpha\int_0^s b(u)\,\dd u\right\}
               \,\dd s \right)^2\\
     &\leq\left(\int_0^t\frac{b(s)^2}{\sigma(s)^2}
                \exp\left\{2\alpha\int_0^s b(u)\,\dd u\right\}\,\dd s\right)
         \!\left(\int_0^t\!\sigma(s)^2
                \exp\left\{-2\alpha\int_0^s b(u)\,\dd u\right\}\!\dd s
          \right).
 \end{align*}
Using \eqref{Fisher2} and that \ $\lim_{t\uparrow T}\int_0^t\vert b(s)\vert\,\dd s=\infty$,
 \ we have \eqref{SEGED_LUSCHGY1}.

\noindent We note that in case of \ $T=\infty$ \ and \ $\sigma\equiv 1$,
 \ Theorem \ref{PRO_CAUCHY} with condition \eqref{Fisher1} replaced by
 \ $\lim_{t\uparrow\infty}\int_0^t\vert b(s)\vert\,\dd s=\infty$ \ was proved by Luschgy
 \cite[Section 4.2]{Lus1}, and hence in Theorem \ref{PRO_CAUCHY} we weaken and generalize
 Luschgy's above mentioned result.
Luschgy's original proof is based on his general theorem (see Theorem 1 in \cite{Lus1}),
 and we note that the conditions of this general theorem are not easy to verify.
Our proof can be considered as a direct one based on limit theorems for local martingales
 (see Theorem \ref{THM_MCLT_RAO2}).
\end{Rem}

The next corollary states that, under the conditions of Theorem \ref{PRO_CAUCHY},
 the MLE of \ $\alpha$ \ is asymptotically normal with an appropriate random normalizing factor.

\begin{Cor}\label{COR_CAUCHY}
Suppose that \ $\alpha \in \RR$ \ such that conditions \eqref{Fisher1} and \eqref{Fisher2}
 are satisfied.
Then
  \[
    \left(\int_0^t\frac{b(u)^2\,(X_u^{(\alpha)})^2}{\sigma(u)^2}\,\dd u\right)^{\frac{1}{2}}
        \left(\widehat\alpha_t^{(X^{(\alpha)})}-\alpha\right)
        \distr\cN(0,1)
        \qquad \text{as\quad $t\uparrow T$.}
  \]
\end{Cor}

\noindent{\bf Proof.}
By \eqref{SEGED20_uj}, we have for all \ $t\in[t_0,T)$,
 \begin{align*}
    \left(\int_0^t\frac{b(u)^2\,(X_u^{(\alpha)})^2}{\sigma(u)^2}\,\dd u\right)^{\frac{1}{2}}
        (\widehat\alpha_t^{(X^{(\alpha)})}-\alpha)
      =\frac{\int_0^t\frac{b(u)X_u^{(\alpha)}}{\sigma(u)}\,\dd B_u}
       {\left(\int_0^t\frac{b(u)^2\,(X_u^{(\alpha)})^2}{\sigma(u)^2}\,\dd u\right)^{\frac{1}{2}}}
      =\frac{\frac{1}{\sqrt{I_{X^{(\alpha)}}(t)}}
             \int_0^t\frac{b(u)X_u^{(\alpha)}}{\sigma(u)}\,\dd B_u}
        {\left(\frac{1}{I_{X^{(\alpha)}}(t)}
          \int_0^t\frac{b(u)^2\,(X_u^{(\alpha)})^2}{\sigma(u)^2}\,\dd u\right)^{\frac{1}{2}}}.
 \end{align*}
\ By the proof of Theorem \ref{PRO_CAUCHY}, we have
 \[
   \PP\left(\lim_{t\uparrow T}
           \frac{1}{I_{X^{(\alpha)}}(t)}\int_0^t\frac{b(u)^2(X_u^{(\alpha)})^2}{\sigma(u)^2}\,\dd u
           =\xi^2\right)=1,
 \]
 where \ $\xi\distre\cN(0,1)$.
\ Hence Theorem \ref{THM_MCLT_RAO2} yields that
 \[
   \left( \frac{\int_0^t \frac{b(s) \, X_s^{(\alpha)}}{\sigma(s)} \, \dd B_s}
               {\sqrt{I_{X^{(\alpha)}}(t)}}, \,
          \frac{\int_0^t
                 \frac{b(s)^2 \, ( X_s^{(\alpha)} )^2}{\sigma(s)^2} \, \dd s}
               {I_{X^{(\alpha)}}(t)} \right)
   \distr \big( |\xi| \eta, \, \xi^2 \big) \qquad
   \text{as \ $t \uparrow T$,}
 \]
 where \ $\eta$ \ is a standard normally distributed random variable
 independent of \ $\vert\xi\vert$.
\ Then, by the continuous mapping theorem,
 we have
 \begin{align*}
   \frac{\frac{1}{\sqrt{I_{X^{(\alpha)}}(t)}}
             \int_0^t\frac{b(u)X_u^{(\alpha)}}{\sigma(u)}\,\dd B_u}
    {\left(\frac{1}{I_{X^{(\alpha)}}(t)}
       \int_0^t\frac{b(u)^2\,(X_u^{(\alpha)})^2}{\sigma(u)^2}\,\dd u\right)^{\frac{1}{2}}}
    \distr \frac{\vert\xi\vert\eta}{\sqrt{\xi^2}}=\eta
     \qquad \text{as \ $t\uparrow T$,}
 \end{align*}
 as stated.
\proofend

Finally we formulate conditions for asymptotic normality.
In the case of \ $T=\infty$ \ and \ $\sigma\equiv1$, \ the corresponding assertion
 has already been formulated and proved in Luschgy \cite[Section 4.2]{Lus1}.
The proof of our more general theorem is the same as in Luschgy \cite[Section 4.2]{Lus1}.

\begin{Thm}\label{PRO_NORMAL}
Suppose that \ $\alpha \in \RR$ \ such that
 \begin{gather}
   \lim_{t \uparrow T} I_{X^{(\alpha)}}(t) = \infty , \label{Fisher3} \\
   \lim_{t \uparrow T}\frac{1}{\sqrt{I_{X^{(\alpha)}}(t)}}
   \frac{b(t)}{\sigma(t)^2}
   \int_0^t\sigma(s)^2 \exp\left\{ 2 \alpha \int_s^t b(v) \, \dd v \right\} \, \dd s=0. \label{Fisher4_new}
  \end{gather}
Then
 \[
   \sqrt{I_{X^{(\alpha)}}(t)} \, ( \halpha_t^{(X^{(\alpha)})} - \alpha )
   \distr \cN(0,1) \qquad
   \text{as \ $t \uparrow T$.}
 \]
\end{Thm}

The next corollary states that, under the conditions of Theorem \ref{PRO_NORMAL},
 the MLE of \ $\alpha$ \ is also asymptotically normal with an appropriate random normalizing factor.

\begin{Cor}\label{COR_NORMAL}
Suppose that \ $\alpha \in \RR$ \ such that conditions \eqref{Fisher3} and \eqref{Fisher4_new}
 are satisfied.
Then
  \[
    \left(\int_0^t\frac{b(u)^2\,(X_u^{(\alpha)})^2}{\sigma(u)^2}\,\dd u\right)^{\frac{1}{2}}
        (\widehat\alpha_t^{(X^{(\alpha)})}-\alpha)
        \distr\cN(0,1)
        \qquad \text{as\quad $t\uparrow T$.}
  \]
\end{Cor}

\noindent{\bf Proof.}
For all \ $t\in[t_0,T)$, \ we have
 \begin{align}\label{SEGED49}
   \begin{split}
     &\left(\int_0^t\frac{b(u)^2\,(X_u^{(\alpha)})^2}{\sigma(u)^2}\,\dd u\right)^{\frac{1}{2}}
        (\widehat\alpha_t^{(X^{(\alpha)})}-\alpha)\\
        &\qquad\qquad
           =\sqrt{I_{X^{(\alpha)}}(t)}(\widehat\alpha_t^{(X^{(\alpha)})}-\alpha)
          \left(\frac{1}{I_{X^{(\alpha)}}(t)}
               \int_0^t\frac{b(u)^2\,(X_u^{(\alpha)})^2}{\sigma(u)^2}\,\dd u\right)^{\frac{1}{2}}.
   \end{split}
 \end{align}
By Theorem \ref{PRO_NORMAL}, we have
 \[
   \sqrt{I_{X^{(\alpha)}}(t)}(\widehat\alpha_t^{(X^{(\alpha)})}-\alpha)
      \distr\cN(0,1)
     \qquad \text{as \ $t\uparrow T$.}
 \]
Moreover, one can show (see, Luschgy \cite[Section 4.2]{Lus1})
 \[
   \frac{1}{I_{X^{(\alpha)}}(t)}\int_0^t\frac{b(u)^2(X_u^{(\alpha)})^2}{\sigma(u)^2}\,\dd u \stoch 1
        \qquad \text{as \ $t\uparrow T$.}
 \]
Note that if \ $\xi_n$, \ $n\in\NN$, \ are nonnegative random variables such that
 \ $\xi_n\stoch 1$ \ as \ $n\to\infty$, \ then \ $\xi_n\distr 1$ \ as \ $n\to\infty$ \ and
 hence \ $\sqrt{\xi_n}\distr 1$ \ as \ $n\to\infty$.
\ Since the limit $1$ is non-random, we have \ $\sqrt{\xi_n}\stoch 1$ \ as \ $n\to\infty$.
\ Hence
 \[
   \left(\frac{1}{I_{X^{(\alpha)}}(t)}\int_0^t\frac{b(u)^2(X_u^{(\alpha)})^2}{\sigma(u)^2}\,\dd u\right)^{1/2}
       \stoch 1 \qquad \text{as \ $t\uparrow T$.}
 \]
By \eqref{SEGED49} and Slutsky's lemma (see, e.g., Lemma 2.8, (ii) in van der Vaart \cite{Vaart}),
 we have the assertion.
\proofend

\begin{Rem}
We note that the sets of those parameters of \ $\alpha$ \ for which  \eqref{Fisher5} and \eqref{Fisher6},
 for which \eqref{Fisher1} and \eqref{Fisher2}, and for which \eqref{Fisher3} and \eqref{Fisher4_new}
 hold are pairwise disjoint.
This is an immediate consequence of the fact that under the conditions of Theorems
 \ref{sing_special}, \ref{PRO_CAUCHY} and \ref{PRO_NORMAL} the asymptotic distributions of
 the MLE of \ $\alpha$ \ are different from each other.
That is, if we can apply one of the Theorems \ref{sing_special}, \ref{PRO_CAUCHY} and \ref{PRO_NORMAL},
 then it is sure that the other two can not be applied.

\noindent We also remark that in general the set of those parameters of \ $\alpha$ \ for which one of the
 Theorems \ref{sing_special}, \ref{PRO_CAUCHY} and \ref{PRO_NORMAL} can be applied is not necessarily
 the whole \ $\RR$.
\ Due to Luschgy \cite[Section 4.2]{Lus1}, if \ $T=\infty$, \ $b(t):=-\ee^{-t}$, \ $t\geq 0$, \ and
 \ $\sigma\equiv1$, \ then \ $\lim_{t\uparrow T}I_{X^{(\alpha)}}(t)=\infty$ \ is not satisfied and
 hence none of the Theorems \ref{sing_special}, \ref{PRO_CAUCHY} and \ref{PRO_NORMAL} can be applied.
\end{Rem}

\section{Consistency}

First we recall a strong law of large numbers which can be applied to
 stochastic integrals.
The following theorem is a modification of Theorem 3.4.6 in Karatzas and Shreve
 \cite{KarShr} (due to Dambis, Dubins and Schwartz),
 see also Theorem 1.6 in Chapter V in Revuz and Yor \cite{RevYor}.
In fact, our next Theorem \ref{DDS} is Exercise 1.18 in Chapter V
 in Revuz and Yor \cite{RevYor}.

\begin{Thm}\label{DDS}
Let \ $T\in(0,\infty]$ \ be fixed and let
\ $\big(\Omega, \cF, (\cF_t)_{t\in[0,T)}, \PP\big)$ \ be a filtered
 probability space satisfying the usual conditions.
Let \ $(M_t)_{t \in [0,T)}$ \ be a continuous local martingale
with respect to the filtration \ $(\cF_t)_{t\in[0,T)}$ \
such that \ $\PP(M_0=0)=1$ \ and
 \ $\PP(\lim_{t \uparrow T} \langle M \rangle_t=\infty)=1$.
\ For each \ $s \in [0,\infty)$, \ define the stopping time
 \[
   \tau_s := \inf \{ t \in [0,T) : \langle M \rangle_t > s \} .
 \]
Then the time-changed process
 \[
   \left( B_s:=M_{\tau_s}, \, \cF_{\tau_s} \right)_{s \geq 0}
 \]
 is a standard Wiener process. In particular, the filtration \ $(\cF_{\tau_s})_{s\geq 0}$ \
 satisfies the usual conditions and
 \[
  \PP\big(M_t=B_{\langle M\rangle_t}\quad\text{for all \ $t\in[0,T)$}\big)=1.
 \]
\end{Thm}

Now we formulate a strong law of large numbers for continuous local
 martingales.
Compare with L\'epingle \cite[Theoreme 1]{Lep} or with $3^\circ)$ in Exercise 1.16
 in Chapter V in Revuz and Yor \cite{RevYor}.
We note that the above mentioned citations are about continuous local
 martingales with time interval \ $[0,\infty)$, \ but they are also valid for
 continuous local martingales with time interval \ $[0,T)$, \ $T\in(0,\infty)$,
 \ with appropriate modifications in the conditions, see as follows.

\begin{Thm}\label{MSLLN}
Let \ $T\in(0,\infty]$ \ be fixed and let
 \ $\big(\Omega, \cF, (\cF_t)_{t\in[0,T)}, \PP\big)$ \ be a filtered probability
 space satisfying the usual conditions.
Let \ $(M_t)_{t \in [0,T)}$ \ be a continuous local martingale
with respect to the filtration \ $(\cF_t)_{t\in[0,T)}$ \
such that \ $\PP(M_0=0)=1$ \ and
\ $\PP(\lim_{t \uparrow T} \langle M \rangle_t = \infty)=1$.
\ Let \ $f : [1, \infty) \to (0, \infty)$ \ be an increasing function such that
 \[
   \int_1^\infty \frac{1}{f(x)^2} \, \dd x < \infty .
 \]
Then
 \[
   \PP \left( \lim_{t \uparrow T} \frac{M_t}{f(\langle M \rangle_t)} = 0 \right)
   = 1 .
 \]
\end{Thm}

Theorem \ref{DDS} has the following consequence on stochastic integrals.

\begin{Thm}\label{DDS_stoch_int}
Let \ $T\in(0,\infty]$ \ be fixed and let \ $\big(\Omega, \cF, (\cF_t)_{t\in[0,T)}, \PP\big)$ \ be
 a filtered probability space satisfying the usual conditions.
Let \ $(M_t)_{t \in [0,T)}$ \ be a continuous local martingale
 with respect to the filtration \ $(\cF_t)_{t\in[0,T)}$ \ such that
 \ $\PP(M_0=0)=1$.
\ Let \ $(\xi_t)_{t\in[0,T)}$ \ be a progressively measurable process such that
 \[
  \PP\left(\int_0^t(\xi_u)^2\,\dd\langle M\rangle_u<\infty\right)=1,\quad t\in[0,T),
 \]
 and
 \begin{align}\label{SEGED_STRONG_CONSISTENCY2}
   \PP\left(\lim_{t\uparrow T}\int_0^t(\xi_u)^2\,\dd\langle M\rangle_u=\infty\right)=1.
 \end{align}
Let
 \[
  \tau_s:=\inf\left\{t\in[0,T):\int_0^t(\xi_u)^2\,\dd\langle M\rangle_u>s\right\},
    \qquad s\geq 0.
 \]
Then the process \ $(\eta_s,\cF_{\tau_s})_{s\geq 0}$, \ defined by
 \[
   \eta_s:=\int_0^{\tau_s}\xi_u\,\dd M_u,\qquad s\geq0,
 \]
 is a standard Wiener process, and
 \begin{align}\label{SEGED_STOCH_INT_SLLN}
   \PP\left(\lim_{t\uparrow T}\frac{\int_0^t\xi_u\,\dd M_u}
         {\int_0^t(\xi_u)^2\,\dd \langle M\rangle_u}=0\right)=1.
 \end{align}
In case of \ $M_t=B_t$, \ $t\in[0,T)$, \ where \ $(B_t)_{t\in[0,T)}$ \ is a standard Wiener
 process, the progressive measurability of \ $(\xi_t)_{t\in[0,T)}$ \ can be relaxed to measurability
 and adaptedness to the filtration \ $(\cF_t)_{t\in[0,T)}$.
\end{Thm}

For historical fidelity, we note that if \ $T=\infty$ \ and \ $M$ \ is a standard Wiener
 process, then Theorem \ref{DDS_stoch_int} was already formulated and proved in Lemma 17.4
 in Liptser and Shiryaev \cite{LipShiII}.
Our proof differs from the original proof of Liptser and Shiryaev.

\noindent{\bf Proof of Theorem \ref{DDS_stoch_int}.}
By Proposition 3.2.24 in Karatzas and Shreve \cite{KarShr}, the process
 $\int_0^t\xi_u\,\dd M_u$, \ $t\in[0,T)$, \ is a continuous, local martingale with respect to
 the filtration \ $(\cF_t)_{t\in[0,T)}$.
\ Since \ $\int_0^t\xi_u\,\dd M_u$, \ $t\in[0,T)$, \ is continuous almost surely, it is
 square integrable.
Moreover, by page 147 in Karatzas and Shreve \cite{KarShr},
 the quadratic variation process of \ $\int_0^t\xi_u\,\dd M_u$, \ $t\in[0,T)$, \ is
 \[
   \int_0^t(\xi_u)^2\,\dd\langle M\rangle_u,\qquad t\in[0,T).
 \]
Hence Theorem \ref{DDS} implies that \ $(\eta_s,\cF_{\tau_s})_{s\geq 0}$ \ is a standard
 Wiener process.
Using condition \eqref{SEGED_STRONG_CONSISTENCY2}, Theorem \ref{MSLLN} implies
 \eqref{SEGED_STOCH_INT_SLLN}.

\noindent In case of \ $M_t=B_t$, \ $t\in[0,T)$, \ Remark 3.2.11 in Karatzas and Shreve \cite{KarShr}
 gives us that the progressively measurability of \ $(\xi_t)_{t\in[0,T)}$ \ can be relaxed to measurability
 and adaptedness to the filtration \ $(\cF_t)_{t\in[0,T)}$. \
\proofend

\begin{Thm}\label{PRO_strongly_consistent2}
Suppose that \ $\alpha\in\RR$ \ such that
\begin{align}\label{SEGED_STRONG_CONSISTENCY}
  \PP\left(\lim_{t\uparrow T}\int_0^t\frac{b(s)^2(X_s^{(\alpha)})^2}{\sigma(s)^2}\,\dd s
       =\infty\right)=1 .
 \end{align}
Then the maximum likelihood estimator \ $\widehat\alpha_t^{(X^{(\alpha)})}$ \ of
 \ $\alpha$ \ is strongly consistent, i.e.,
 \[
  \PP\big(\lim_{t\uparrow T}\widehat\alpha_t^{(X^{(\alpha)})}=\alpha\big)=1.
 \]
\end{Thm}

\noindent{\bf Proof.}
Using \eqref{SEGED20_uj} and \eqref{SEGED_STRONG_CONSISTENCY}, Theorem \ref{DDS_stoch_int}
 yields the assertion.
\proofend

Note that in the case of an Ornstein-Uhlenbeck process, condition \eqref{SEGED_STRONG_CONSISTENCY}
 is satisfied for all \ $\alpha\in\RR$ \ (see, e.g., Liptser and Shiryaev \cite[(17.57)]{LipShiII}),
 and hence  in this case the strong consistency of the MLE of \ $\alpha$ \ is an immediate
 consequence of Theorem \ref{PRO_strongly_consistent2}.

We also remark that if the conditions of Theorem \ref{sing_special} or Theorem
 \ref{PRO_CAUCHY} or Theorem \ref{PRO_NORMAL} are satisfied then weak consistency of the
 MLE \ of \ $\alpha$ \ holds.

\section{Perturbation of the drift}\label{Perturbation of the drift}

Let \ $T \in (0, \infty]$ \ be fixed.
Let \ $b: [0,T) \to \RR$ \ and \ $\sigma:[0,T) \to \RR$ \ be continuous
 functions.
Suppose that \ $\sigma(t) > 0$ \ for all \ $t\in[0,T)$, \ and there exists
 \ $t_0 \in (0,T)$ \ such that \ $b(t)\ne 0$ \ for all \ $t\in[t_0,T)$.
\ Let \ $a:\RR\to\RR$ \ be a function such that
 \[
   a(x) = x + r(x) , \qquad x \in \RR ,
 \]
 where
 \[
   |r(x)| \leq L (1 + |x|^\gamma ) , \qquad x \in \RR ,
 \]
 with some \ $L \geq 0$ \ and \ $\gamma \in [0,1)$, \ and \ $r$ \ satisfies
 the global Lipschitz condition
 \begin{align}\label{CONDITION_ON_r}
   |r(x) - r(y)| \leq M |x-y| , \qquad x,y \in \RR ,
 \end{align}
 with some \ $M\geq0$.
\ Note that continuity of \ $r$ \ implies continuity of \ $a$.
\ For all \ $\alpha \in \RR$, \ let us consider the SDE \eqref{perturbed_SDE}.
Note that the drift and diffusion coefficients of the SDE \eqref{perturbed_SDE}
 satisfy the local Lipschitz condition and the linear growth condition
 (see, e.g., Jacod and Shiryaev \cite[Theorem 2.32, Chapter III]{JacShi}).
Again by Jacod and Shiryaev \cite[Theorem 2.32, Chapter III]{JacShi}, the SDE
 \eqref{perturbed_SDE} has a unique strong solution.
Note also that \ $(Y_t^{(\alpha)})_{t \in [0,T)}$ \ has continuous sample paths by the
 definition of strong solution, see, e.g., Jacod and Shiryaev
 \cite[Definition 2.24, Chapter III]{JacShi}.
For all \ $\alpha \in \RR$ \ and \ $t \in (0,T)$, \ let \ $\PP_{Y^{(\alpha)},\,t}$
 \ denote the distribution of the process \ $(Y_s^{(\alpha)})_{s \in [0,\,t]}$ \ on
 \ $\big(C([0,t]),\cB(C([0,t]))\big)$.
\ The measures \ $\PP_{Y^{(\alpha)},\,t}$ \ and \ $\PP_{Y^{(0)},\,t}$ \ are equivalent and
 \[
   \frac{\dd \PP_{Y^{(\alpha)},\,t}}{\dd \PP_{Y^{(0)},\,t}}
    \left(Y^{(\alpha)}\big\vert_{[0,t]}\right)
   = \exp \left\{ \alpha
                   \int_0^t
                    \frac{b(s)}{\sigma(s)^2} \, a\big(Y_s^{(\alpha)}\big) \,
                     \dd Y_s^{(\alpha)}
                   -\frac{\alpha^2}{2}
                    \int_0^t
                     \frac{b(s)^2}{\sigma(s)^2} \, a\big( Y_s^{(\alpha)} \big)^2
                      \, \dd s \right\} ,
 \]
 see Liptser and Shiryaev \cite[Theorem 7.20]{LipShiI}.

The MLE \ $\widehat\alpha_t^{(Y^{(\alpha)})}$ \ of \ $\alpha$ \ based on the
 observation \ $(Y_s^{(\alpha)})_{s\in[0,\,t]}$ \ is defined by
 \[
    \widehat\alpha_t^{(Y^{(\alpha)})}
      :=\argmax_{\alpha\in\RR} \,
         \ln\left(\frac{\dd\,\PP_{Y^{(\alpha)},\,t}}{\dd\,\PP_{Y^{(0)},\,t}}
                 \left(Y^{(\alpha)}\big\vert_{[0,t]}\right)\right).
 \]
If \ $\omega\in\Omega$ \ such that
 \ $\int_0^t\frac{b(s)^2}{\sigma(s)^2}a(Y_s^{(\alpha)}(\omega))^2\,\dd s=0$
 \ and \ $\left(\int_0^t\frac{b(s)}{\sigma(s)^2} \, a\big(Y_s^{(\alpha)}\big) \,
        \dd Y_s^{(\alpha)}\right)(\omega)\ne0$, \ then
 \[
   \sup_{\alpha\in\RR}\,
    \ln\left(\frac{\dd\,\PP_{Y^{(\alpha)},\,t}}{\dd\,\PP_{Y^{(0)},\,t}}
             \left(Y^{(\alpha)}\big\vert_{[0,t]}\right)\right)(\omega)
   =\infty ,
 \]
 which yields that \ $\widehat\alpha_t^{(Y^{(\alpha)})}(\omega)$ \
 does not exist.
If condition
 \begin{align}\label{CON_MLE}
    \PP\left(\lim_{t\uparrow T}
           \int_0^t\frac{b(s)^2}{\sigma(s)^2}a(Y^{(\alpha)}_s)^2\,\dd s>0
           \right)=1
 \end{align}
 holds, then the MLE \ $\widehat\alpha_t^{(Y^{(\alpha)})}$ \ of \ $\alpha$ \ based on
 the observation \ $(Y_s^{(\alpha)})_{s\in[0,\,t]}$ \ exists asymptotically as \ $t\uparrow T$
 \ with probability one.
(Note that in case of \ $r\equiv0$ \ condition \eqref{SEGED_STRONG_CONSISTENCY} yields \eqref{CON_MLE}.)
Hence
 \[
   \widehat\alpha_t^{(Y^{(\alpha)})}
      =\frac{\int_0^t\frac{b(s)a(Y_s^{(\alpha)})}{\sigma(s)^2}
              \,\dd Y_s^{(\alpha)}}
             {\int_0^t\frac{b(s)^2a(Y_s^{(\alpha)})^2}{\sigma(s)^2}\,\dd s}
 \]
 holds asymptotically as \ $t\uparrow T$ \ with probability one.
To be more precise, if condition \eqref{CON_MLE} holds, there exists an event
 \ $A\in\cF$ \ such that \ $\PP(A)=1$ \ and for all \ $\omega\in A$ \ there
 exists a \ $t(\omega)\in[0,T)$ \ with the property that
 \ $\widehat\alpha_t^{(Y^{(\alpha)})}(\omega)$ \ exists for all
 \ $t\in[t(\omega),T)$ \ and
 \[
   \widehat\alpha_t^{(Y^{(\alpha)})}(\omega)
      =\frac{\left(\int_0^t\frac{b(s)a(Y_s^{(\alpha)})}{\sigma(s)^2}
              \,\dd Y_s^{(\alpha)}\right)(\omega)}
             {\int_0^t\frac{b(s)^2a(Y_s^{(\alpha)}(\omega))^2}{\sigma(s)^2}\,\dd s}.
 \]
In all what follows, by the expression `exists/holds asymptotically as
 \ $t\uparrow T$ \ with probability one' we mean the above property.
Using the SDE \eqref{perturbed_SDE}, we have for all \ $\alpha\in\RR$,
 \begin{align*}
   \widehat\alpha_t^{(Y^{(\alpha)})}-\alpha
     =\frac{\int_0^t\frac{b(s)a(Y_s^{(\alpha)})}{\sigma(s)}
              \,\dd B_s}
             {\int_0^t\frac{b(s)^2a(Y_s^{(\alpha)})^2}{\sigma(s)^2}\,\dd s}
 \end{align*}
 holds asymptotically as \ $t\uparrow T$ \ with probability one.

The following lemma gives a sufficient condition under which \eqref{CON_MLE} is satisfied
 for all \ $\alpha\in\RR$.

\begin{Lem}\label{LEMMA_PERTURB_FELTETEL}
If \ $\lim_{t\uparrow T}\int_0^t\sigma(s)^2\,\dd s=\infty$, \ then
 \eqref{CON_MLE} is satisfied for all \ $\alpha\in\RR$.
\end{Lem}

\noindent{\bf Proof.}
We follow the ideas of the proof of Lemma 3.2 in Dietz and Kutoyants \cite{DieKut}.
Let \ $\alpha\in\RR$ \ be fixed.
On the contrary, let us suppose that \ $\PP(A_1)>0$, \ where
 \[
    A_1:=\left\{\omega\in\Omega :
                 \lim_{t\uparrow T}\int_0^t\frac{b(s)^2}{\sigma(s)^2}\,
                    a\big(Y_s^{(\alpha)}(\omega)\big)^2\,\dd s
        =0\right\}.
 \]
Then for all \ $t\in[0,T)$ \ and \ $\omega\in A_1$, \ we have
 \[
     \int_0^t\frac{b(s)^2}{\sigma(s)^2}\,a\big(Y_s^{(\alpha)}(\omega)\big)^2\,\dd s=0.
 \]
Since \ $b$, \ $\sigma$, \ $a$ \ are continuous and
 \ $Y^{(\alpha)}_.(\omega)$ \ is also continuous on \ $[0,T)$ \ for all \ $\omega\in\Omega$,
 \ we have
 \[
   b(t)a\big(Y_t^{(\alpha)}(\omega)\big)=0,
      \qquad \forall\;\,t\in[0,T),\quad \forall\;\,\omega\in A_1.
 \]
This yields that \ $A_1\subset A_2$, \ where
 \[
    A_2:=\left\{\omega\in\Omega :
                 b(t)a\big(Y_t^{(\alpha)}(\omega)\big)=0,
      \quad \forall\;\,t\in[0,T)\right\}.
 \]
Let \ $Z:=\{x\in\RR:a(x)=0\}$.
\ We show that \ $Z$ \ is compact.
First we check that \ $\lim_{x\to\pm\infty}a(x)=\pm\infty$.
\ Since
 \[
    \left\vert\frac{r(x)}{x}\right\vert\leq L\left(\frac{1}{\vert x\vert}+\vert x\vert^{\gamma-1}\right)
                              \to 0
       \qquad \text{as} \quad x\to\pm\infty,
 \]
 we have
 \[
    \lim_{x\to\pm\infty}a(x)=\lim_{x\to\pm\infty}x\left(1+\frac{r(x)}{x}\right)=\pm\infty.
 \]
Hence, using also that \ $a$ \ is continuous, we have \ $Z$ \ is compact.
Using that \ $b(t)\ne0$ \ for all \ $t\in[t_0,T)$, \ we have
 \[
    a\big(Y_t^{(\alpha)}(\omega)\big)=0,
      \qquad \forall\;\,\omega\in A_2,\quad \forall\;\,t\in[t_0,T),
 \]
 i.e., \ $Y_t^{(\alpha)}(\omega)\in Z$ \ for all \ $\omega\in A_2$ \ and for all \ $t\in[t_0,T)$.
\ By the SDE \eqref{perturbed_SDE}, we have
 \[
   Y_t^{(\alpha)}(\omega)=\left(\int_0^t\sigma(s)\,\dd B_s\right)(\omega),
           \qquad \forall\;\omega\in A_2,\quad\forall\;t\in[0,T),
 \]
 and hence
 \[
   a\left(\left(\int_0^t\sigma(s)\,\dd B_s\right)(\omega)\right)=0,
         \qquad \forall\;\omega\in A_2,\quad\forall\;t\in[t_0,T),
 \]
 i.e., \ $\left(\int_0^t\sigma(s)\,\dd B_s\right)(\omega)\in Z$
 \ for all \ $\omega\in A_2$ and for all \ $t\in[t_0,T)$.
\ Then
 \[
    0<\PP(A_1)\leq \PP(A_2)
      \leq \PP\left(\left\{\omega\in \Omega:\left(\int_0^t\sigma(s)\,\dd B_s\right)(\omega)\in Z,
                       \quad\forall\;t\in[t_0,T)\right\}\right).
 \]
This leads us to a contradiction. Indeed, the Gauss process
 \ $\left(\int_0^t\sigma(s)\,\dd B_s\right)_{t\in[0,T)}$ \ has expectation
 function \ $0$ \ and variance function \ $\int_0^t\sigma(s)^2\,\dd s$, \ $t\in[0,T)$.
\ Using that \ $Z$ \ is compact, there exists \ $K>0$ \ such that \ $\vert x\vert<K$ \ for
 all \ $x\in Z$.
\ Hence
 \begin{align}\label{SEGED_PERTURB2}
  \begin{split}
   0&<\PP\left(\left\{\omega\in\Omega:
             \left\vert \left(\int_0^t\sigma(s)\,\dd B_s\right)(\omega)\right\vert<K,
                \;\; \forall\;\,t\in[t_0,T)\right\}\right)\\
    &\leq \PP\left(\left\vert\int_0^t\sigma(s)\,\dd B_s\right\vert<K\right),
          \quad \forall\;\,t\in[t_0,T).
   \end{split}
 \end{align}
Using that, by our assumption, \ $\lim_{t\uparrow T}\int_0^t\sigma(s)^2\,\dd s=\infty$ \ and that
 \[
   \int_0^t\sigma(s)\,\dd B_s \distre \cN\left(0,\int_0^t\sigma(s)^2\,\dd s\right),
    \qquad t\in[0,T),
 \]
 we get
 \[
  \frac{\int_0^t\sigma(s)\,\dd B_s}{\sqrt{\int_0^t\sigma(s)^2\,\dd s}}
     \distre\cN(0,1),
     \qquad t\in(0,T).
 \]
Hence
 \begin{align*}
    \lim_{t\uparrow T}
      \PP\left(\left\vert\int_0^t\sigma(s)\,\dd B_s\right\vert<K\right)
       =\lim_{t\uparrow T}
             \PP\left(\frac{\left\vert\int_0^t\sigma(s)\,\dd B_s\right\vert}{\sqrt{\int_0^t\sigma(s)^2\,\dd s}}
                <\frac{K}{\sqrt{\int_0^t\sigma(s)^2\,\dd s}}\right)
       =\PP(\vert\xi\vert<0)=0,
 \end{align*}
 where \ $\xi$ \ is a standard normally distributed random variable.
Here the last but one equality follows by the fact that if \ $F_n$, \ $n\in\NN$, \
 are distribution functions such that \ $\lim_{n\to\infty}F_n(x)=F(x)$ \
 for all \ $x\in\RR$, \ where \ $F$ \ is a continuous distribution function,
 then for all sequences \ $(x_n)_{n\in\NN}$ \ for which \ $\lim_{n\to\infty}x_n=x\in\RR$,
 \ we have \ $\lim_{n\to\infty}F_n(x_n)=F(x)$.
\ By \eqref{SEGED_PERTURB2}, we arrive at a contradiction.
\proofend

In the next remark we give an example for \ $\alpha$, \ $b$, \ $r$ \ and \ $\sigma$ \
 for which condition \eqref{CON_MLE} does not hold, and also give an example
 for which it holds.

\begin{Rem}
We will give an example for \ $\alpha$, \ $b$, \ $r$ \ and \ $\sigma$ \ such that
 for all \ $t\in[0,T)$,
 \[
   \PP\left(\int_0^t\frac{b(s)^2}{\sigma(s)^2}a(Y^{(\alpha)}_s)^2\,\dd s=0
           \right)>0.
 \]
In this case, for all \ $t\in(0,T)$, \ the MLE \ $\widehat\alpha_t^{(Y^{(\alpha)})}$
 \ of \ $\alpha$ \ exists only with probability less than one.
We note that in our example condition \eqref{CON_MLE} will not hold, and hence the MLE
 of \ $\alpha$ \ will exist asymptotically as \ $t\uparrow T$ \ only with probability less than one.
We also give an example for \ $\alpha$, \ $b$, \ $r$ \ and \ $\sigma$ \ such that for all
 \ $t\in(0,T)$,
 \[
   \PP\left(\int_0^t\frac{b(s)^2}{\sigma(s)^2}a(Y^{(\alpha)}_s)^2\,\dd s=0
           \right)=0.
 \]
In this case, for all \ $t\in(0,T)$, \ the MLE \ $\widehat\alpha_t^{(Y^{(\alpha)})}$
 \ exists with probability one, and condition \eqref{CON_MLE} holds trivially.

First we consider the case \ $T\in(0,\infty)$.
\ Let \ $b(t):=1$, \ $t\in[0,T)$, \ $\sigma(t):=\frac{1}{\sqrt{T-t}}$, \ $t\in[0,T)$, \ and
 \ $r(t):=0$, \ $t\in[0,T)$.
\ Since \begin{align*}
   \lim_{t\uparrow T}\int_0^t\sigma(s)^2\,\dd s
       =\lim_{t\uparrow T}\int_0^t\frac{1}{T-s}\,\dd s
       =\infty,
 \end{align*}
 by Lemma \ref{LEMMA_PERTURB_FELTETEL}, we get condition \eqref{CON_MLE} is satisfied
 for all \ $\alpha\in\RR$.

\noindent In what follows we give an example for \ $\alpha$, \ $b$, \ $r$ \ and \ $\sigma$
 \ such that for all \ $t\in[0,T)$,
 \[
   \PP\left(\int_0^t\frac{b(s)^2}{\sigma(s)^2}a(Y^{(\alpha)}_s)^2\,\dd s=0\right)
     >0.
 \]
In fact, we just reformulate Remark 3.1 in Dietz and Kutoyants \cite{DieKut}
 which is originally stated for the time interval \ $[0,\infty)$.
\ Let \ $b(t):=1$, \ $t\in[0,T)$, \ $\sigma(t):=1$, \ $t\in[0,T)$, \ and
 \[
    r(x):=\begin{cases}
            \frac{1}{1+(x+1)^2} & \text{\quad if \ $x<-1$,}\\
            -x & \text{\quad if \ $-1\leq x<1$,}\\
            -\frac{1}{1+(x-1)^2} & \text{\quad if \ $1\leq x$.}\\
         \end{cases}
 \]
Note that in this case \ $\lim_{t\uparrow T}\int_0^t\sigma(s)^2\,\dd s=T<\infty$,
 \ and hence one can not use Lemma \ref{LEMMA_PERTURB_FELTETEL} for proving \eqref{CON_MLE}.
It will turn out that \eqref{CON_MLE} is not satisfied for \ $\alpha=1$.
\ Clearly, \ $r$ \ is continuous, piecewise continuously differentiable and has everywhere
 left and right derivatives.
Moreover, \ $\vert r(x)\vert\leq 1$, \ $x\in\RR$, \ and all of its (one-sided) derivatives
 are bounded by 1.
Therefore, \ $\vert r(x)\vert\leq L(1+\vert x\vert^\gamma)$, \ $x\in\RR$, \ and
 \ $\vert r(x)-r(y)\vert\leq M\vert x-y\vert,$ $x,y\in\RR$, \ with \ $L:=\frac{1}{2}$,
 \ $\gamma:=0$ \ and \ $M:=1$.
\ (The fact that one can choose \ $M$ \ to be 1 follows from Lagrange's theorem.
Note that \ $r$ \ is not differentiable everywhere, but we can apply Lagrange's theorem
 on different subintervals of \ $\RR$ \ separately, where \ $r$ \ is differentiable.)
Let \ $\alpha:=1$.
\ Then, by the SDE \eqref{perturbed_SDE},
 \begin{align}\label{perturbed_SDE_SOL_SPEC}
    Y_t^{(1)}=\int_0^ta(Y_s^{(1)})\,\dd s+B_t,\qquad t\in[0,T).
 \end{align}
Let us define the random variable \ $\tau$ \ by
 \[
   \tau(\omega):=
      \begin{cases}
         \inf\big\{t\in[0,T):\vert Y_t^{(1)}(\omega)\vert\geq 1\big\}
              & \text{if \ $\exists$ \ $t\in[0,T):$ $\vert Y_t^{(1)}(\omega)\vert\geq 1,$}\\
         T  & \text{if \ $\vert Y_t^{(1)}(\omega)\vert<1,$ $\forall$ $t\in[0,T).$}
       \end{cases}
 \]
Since \ $Y_0^{(1)}=0$ \ and \ $(Y_t^{(1)}(\omega))_{t\in[0,T)}$ \ is continuous for all
 \ $\omega\in\Omega$, \ we have \ $\PP(\tau>0)=1$, \ and if \ $\tau(\omega)<T$, \ then
 \ $\vert Y_{\tau(\omega)}^{(1)}(\omega)\vert=1$.
\ By the definition of \ $\tau$, \ we have \ $\vert Y_t^{(1)}(\omega)\vert<1$ \ for all
 \ $0\leq t<\tau(\omega)$.
\ Hence, using that \ $a(x)=x+r(x)=0$, \ $\vert x\vert\leq 1$, \ we have \ $a(Y_t^{(1)}(\omega))=0$
 \ for all \ $0\leq t<\tau(\omega)$, \ and then
 \[
      \int_0^t a(Y_s^{(1)}(\omega))\,\dd s
       =\int_0^t a(Y_s^{(1)}(\omega))^2\,\dd s=0,
        \qquad 0\leq t<\tau(\omega).
 \]
Hence, by \eqref{perturbed_SDE_SOL_SPEC}, we have
 \ $Y_t^{(1)}(\omega)=B_t(\omega)$, \ $0\leq t<\tau(\omega)$.
\ Note that if \ $\tau(\omega)<T$, \ then we also have
 \ $Y_{\tau(\omega)}^{(1)}(\omega)=B_{\tau(\omega)}(\omega)$ \ and hence
 \ $\vert B_{\tau(\omega)}(\omega)\vert=1$.
\ Let us define the random variable \ $\kappa$ \ by
 \[
  \kappa(\omega):=\inf\Big\{t\in[0,\infty):\vert B_t(\omega)\vert\geq 1\Big\}.
 \]
Hence, if \ $\tau(\omega)<T$, \ we get \ $\kappa(\omega)=\tau(\omega)$, \ and
 if \ $\tau(\omega)=T$, \ then \ $\kappa(\omega)\geq T$.
\ By formula 2.0.2 on page 163 in Borodin and Salminen \cite{BorSal},
 \ $\kappa$ \ is unbounded and \ $\PP(\kappa<\infty)=1$.
\ Consequently,
 \begin{align*}
    0<\PP(\kappa\geq t)
      &=\PP(\{\kappa\geq t\}\cap\{\tau<T\})
        +\PP(\{\kappa\geq t\}\cap\{\tau=T\})\\
      &=\PP(\{\tau\geq t\}\cap\{\tau<T\})
        +\PP(\tau=T)
     \leq 2\PP\left(\int_0^t a(Y_s^{(1)})^2\,\dd s=0\right),
     \quad \forall\;\;t\in[0,T),
 \end{align*}
 as desired.
This also implies that
 \begin{align*}
   0<\PP(\kappa\geq T)
   \leq 2\lim_{t\uparrow T} \PP\left(\int_0^t a(Y_s^{(1)})^2\,\dd s=0\right)
    =2\PP\left(\lim_{t\uparrow T}\int_0^t a(Y_s^{(1)})^2\,\dd s=0\right),
 \end{align*}
 hence \eqref{CON_MLE} is not satisfied for \ $\alpha=1$.

Now we consider the case \ $T=\infty$.
\ Let \ $b(t):=1$, \ $t\geq 0$, \ $\sigma(t):=1$, \ $t\geq 0$, \ and \ $r(t):=0$, \ $t\geq 0$.
\ Since
 \[
  \lim_{t\uparrow \infty}\int_0^t\sigma(s)^2\,\dd s=\lim_{t\uparrow \infty} t=\infty,
 \]
 by Lemma \ref{LEMMA_PERTURB_FELTETEL}, we get \eqref{CON_MLE} holds for all \ $\alpha\in\RR$.
\ Remark 3.1 in Dietz and Kutoyants \cite{DieKut}
 (which we already reformulated for the case \ $T\in(0,\infty)$) \
 gives an example for \ $\alpha$, \ $b$, \ $r$ \ and \ $\sigma$
 \ such that
 \[
   \PP\left(\int_0^t\frac{b(s)^2}{\sigma(s)^2}a(Y^{(\alpha)}_s)^2\,\dd s=0\right)
      >0, \qquad t\in[0,\infty).
 \]
In this example we also have
 \[
   \lim_{t\uparrow \infty}
     \int_0^t\sigma(s)^2\,\dd s
    =\lim_{t\uparrow \infty} t
      =\infty,
 \]
 and hence, by Lemma \ref{LEMMA_PERTURB_FELTETEL}, we have \eqref{CON_MLE} holds for all \ $\alpha\in\RR$.

\noindent In case of \ $T=\infty$ \ we are not able to give an example for \ $\alpha$, \ $b$, \ $r$
 \ and \ $\sigma$ \ such that
 \[
   \PP\left(\int_0^t\frac{b(s)^2}{\sigma(s)^2}a(Y^{(\alpha)}_s)^2\,\dd s=0\right)>0,
     \qquad \forall\,\;t\in[0,\infty),
 \]
 and condition \eqref{CON_MLE} is not satisfied.
For such an example, by Proposition 1.26 in Chapter IV, Proposition 1.8 in Chapter V in Revuz
 and Yor \cite{RevYor} and Lemma \ref{LEMMA_PERTURB_FELTETEL}, it is necessary to have
 \[
   \lim_{t\uparrow \infty}\int_0^t\sigma(s)^2\,\dd s<\infty
   \qquad\text{and}\qquad
   \PP\left(\lim_{t\uparrow \infty}\int_0^t\sigma(s)\,\dd B_s=\zeta\right)=1,
 \]
 where \ $\zeta$ \ is a normally distributed random variable with mean 0 and with variance
 \ $\int_0^{\infty}\sigma(s)^2\,\dd s$.
\end{Rem}

For all \ $t\in(0,T)$, \ the Fisher information for \ $\alpha$ \ contained in
 the observation \ $(Y_s^{(\alpha)})_{s\in[0,\,t]}$ \ is defined by
 \[
   I_{Y^{(\alpha)}}(t)
   := \EE \left( \frac{\partial}{\partial\alpha}
                 \ln \left(\frac{\dd \PP_{Y^{(\alpha)},\,t}}{\dd \PP_{Y^{(0)},\,t}}
                  \left( Y^{(\alpha)} \big \vert_{[0,t]} \right) \right) \right)^2
   = \int_0^t
      \frac{b(s)^2}{\sigma(s)^2} \, \EE a\big( Y_s^{(\alpha)} \big)^2  \, \dd s ,
 \]
 where the last equality follows by the SDE \eqref{perturbed_SDE} and
 Karatzas and Shreve \cite[Proposition 3.2.10]{KarShr}.
Note that \ $I_{Y^{(\alpha)}}(t)\geq 0$ \ for all \ $t\in[0,T)$, \ but in general
\ $I_{Y^{(\alpha)}}(t)>0$, \ $\forall$ $t\in[0,T)$, \ does not hold necessarily.

\begin{Thm}  \label{sing_perturb}
Suppose that \ $\alpha\in\RR$ \ such that
 \begin{gather}
   \lim_{t \uparrow T} I_{X^{(\alpha)}}(t)=\infty,\label{Fisherp5}\\[1mm]
  \lim_{t \uparrow T}
   \frac{b(t)}{\sigma(t)^2}
   \exp\left\{ 2 \alpha \int_0^t b(w) \, \dd w \right\}
      =C\in\RR\setminus\{0\}, \label{Fisherp8}
 \end{gather}
 and \ $\sign(\alpha)=\sign(C)$ \ or \ $\alpha=0$.
Then
 \[
   \sqrt{I_{Y^{(\alpha)}}(t)} \, ( \halpha_t^{(Y^{(\alpha)})} - \alpha )
   \distr \frac{\sign(C)}{\sqrt{2}} \,
          \frac{\int_0^1 W_s \, \dd W_s}{\int_0^1 (W_s)^2 \, \dd s} \qquad
   \text{as \ $t \uparrow T$.}
 \]
\end{Thm}

Note that conditions \eqref{Fisherp5} and \eqref{Fisherp8} do not contain the
 function \ $r$.
\ For the proof of Theorem \ref{sing_perturb}, we need a generalization of
 Gr\"onwall's inequality. Our generalization can be considered as a slight improvement
 of Bainov and Simeonov \cite[Lemma 1.1]{BaiSem}. The proof goes along the same lines.

\begin{Lem}[\textbf{A generalization of Gr\"onwall's inequality}] \label{Gronwall}
Let \ $s_0,s_1\in\RR$ \ with \ $s_0<s_1$, \ let
 \ $\varphi:[s_0,s_1]\to[0,\infty)$ \ and \ $\psi_2:[s_0,s_1]\to[0,\infty)$ \ be
 continuous functions, and let \ $\psi_1:[s_0,s_1]\to\RR$ \ be a continuously
 differentiable function. Suppose that
 \[
   \varphi(s) \leq \psi_1(s) + \int_{s_0}^s \psi_2(u) \varphi(u) \, \dd u ,
   \qquad s \in [s_0,s_1] .
 \]
Then
 \[
   \varphi(s)
   \leq \psi_1(s_0) \exp\left\{ \int_{s_0}^s \psi_2(u) \, \dd u \right\}
        + \int_{s_0}^s
           \psi_1'(u) \exp\left\{ \int_u^s \psi_2(v) \, \dd v \right\} \dd u ,
   \qquad s \in [s_0,s_1] .
 \]
\end{Lem}

\noindent{\bf Proof of Theorem \ref{sing_perturb}.}
Note that condition \eqref{Fisherp8} yields that there exists \ $t_0\in(0,T)$ \ such that
 \ $b(t)\ne0$ \ for all \ $t\in[t_0,T)$.
\ First we check that \ $\lim_{t\uparrow T}\int_0^t \sigma(s)^2 \dd s = \infty$.
\ By \eqref{Fisherp8}, there exist \ $c_1>0$, \ $c_2>0$ \ and \ $t_1\in[t_0,T)$ \ such that
 \eqref{liminfsup} is satisfied.
Hence for all \ $t\in[t_1,T)$,
 \[
  \int_0^t \sigma(s)^2\,\dd s
    \geq \int_0^{t_1} \sigma(s)^2\,\dd s
           +\int_{t_1}^t c_1 \vert b(s) \vert \exp\left\{2\alpha\int_0^s b(w)\,\dd w\right\}\,\dd s.
 \]
By Lemma \ref{LEMMA_EKVIVALENS_SINGULAR}, we have
 \begin{align}\label{Fisherp6}
   \lim_{t\uparrow T}\int_0^t \vert b(s)\vert \, \dd s = \infty.
 \end{align}
If \ $\alpha>0$ \ and \ $C>0$, \ then \ $b(t)>0$ \ for all \ $t\in[t_1,T)$ \ and
 \ $\lim_{t\uparrow T}\int_0^t b(s)\,\dd s=\infty$.
\ Hence
 \[
    \lim_{t\uparrow T}
       \int_{t_1}^t c_1 \vert b(s) \vert \exp\left\{2\alpha\int_0^s b(w)\,\dd w\right\}\,\dd s
       =\infty,
 \]
 which yields \ $\lim_{t\uparrow T}\int_0^t \sigma(s)^2 \dd s = \infty$.

\noindent If \ $\alpha=0$, \  by \eqref{Fisherp6}, we have
 \begin{align*}
   \lim_{t\uparrow T}
      \int_{t_1}^t c_1 \vert b(s)\vert \exp\left\{2\alpha\int_0^s b(w)\,\dd w\right\}\,\dd s
     =c_1 \lim_{t\uparrow T} \int_{t_1}^t\vert b(s)\vert\,\dd s
     =\infty,
 \end{align*}
 which yields \ $\lim_{t\uparrow T}\int_0^t \sigma(s)^2 \dd s = \infty$.

\noindent If \ $\alpha<0$ \ and \ $C<0$, \ then \ $b(t)<0$ \ for all \ $t\in[t_1,T)$ \
 and \ $\lim_{t\uparrow T}\int_0^t b(s)\,\dd s=-\infty$.
\ Hence
 \[
    \lim_{t\uparrow T}
       \int_{t_1}^t c_1 \vert b(s) \vert \exp\left\{2\alpha\int_0^s b(w)\,\dd w\right\}\,\dd s
       =\infty,
 \]
 which yields \ $\lim_{t\uparrow T}\int_0^t \sigma(s)^2 \dd s = \infty$.
\ By Lemma \ref{LEMMA_PERTURB_FELTETEL}, condition \eqref{CON_MLE} holds for
 all \ $\alpha$ \ satisfying the assumptions of the present theorem.
Hence the MLE \ $\widehat\alpha_t^{(Y^{(\alpha)})}$ \ of
 \ $\alpha$ \ based on the observation \ $(Y_s^{(\alpha)})_{s\in[0,\,t]}$ \ exists asymptotically
 as \ $t\uparrow T$ \ with probability one.

Consider the SDEs \eqref{perturbed_SDE} and \eqref{special_SDE}.
Introduce the stochastic process
 \[
   \Delta_t^{(\alpha)} := Y_t^{(\alpha)} - X_t^{(\alpha)} , \qquad t \in [0,T) .
 \]
This process satisfies the ordinary differential equation
 \[
   \left\{
    \begin{aligned}
     \dd \Delta_t^{(\alpha)}
     &= \alpha \, b(t) \Delta_t^{(\alpha)} \, \dd t
        + \alpha b(t) r(Y_t^{(\alpha)}) \, \dd t ,
     \qquad t \in [0,T) , \\
     \Delta_0^{(\alpha)} &= 0 ,
    \end{aligned}
   \right .
 \]
 having the unique solution
 \begin{equation} \label{SolDelta}
   \Delta_t^{(\alpha)}
   = \alpha
     \int_0^t
      r(Y_s^{(\alpha)}) b(s) \exp\left\{ \alpha \int_s^t b(u) \, \dd u \right\}
      \dd s , \qquad t \in [0,T) .
 \end{equation}
Using the decomposition
 \begin{equation} \label{DECO}
   a(Y_t^{(\alpha)})
   = Y_t^{(\alpha)} + r(Y_t^{(\alpha)})
   = X_t^{(\alpha)} + \Delta_t^{(\alpha)} + r(Y_t^{(\alpha)}) , \qquad t \in [0,T) ,
 \end{equation}
 we get
 \[
   \halpha_t^{(Y^{(\alpha)})} - \alpha
   = \frac{\int_0^t \frac{b(s)}{\sigma(s)} \, X_s^{(\alpha)} \, \dd B_s
           + J_t^{(1)} + J_t^{(2)}}
          {\int_0^t \frac{b(s)^2}{\sigma(s)^2} \, ( X_s^{(\alpha)} )^2 \, \dd s
           + J_t^{(3)} + J_t^{(4)} + J_t^{(5)} + J_t^{(6)}}
 \]
 holds asymptotically as \ $t\uparrow T$ \ with probability one, where for all \ $t\in[0,T),$ \
 \begin{align*}
  &J_t^{(1)} := \int_0^t
               \frac{b(s)}{\sigma(s)} \, \Delta_s^{(\alpha)} \, \dd B_s , \qquad
   &J_t^{(2)} := \int_0^t \frac{b(s)}{\sigma(s)} \, r(Y_s^{(\alpha)}) \, \dd B_s , \\[1mm]
  &J_t^{(3)} := 2 \int_0^t
                 \frac{b(s)^2}{\sigma(s)^2} \,
                 X_s^{(\alpha)} \Delta_s^{(\alpha)} \, \dd s , \qquad
   &J_t^{(4)} := \int_0^t
               \frac{b(s)^2}{\sigma(s)^2} \,
               ( \Delta_s^{(\alpha)} )^2 \, \dd s , \\[1mm]
  &J_t^{(5)} := 2 \int_0^t
                 \frac{b(s)^2}{\sigma(s)^2} \,
                 Y_s^{(\alpha)} r(Y_s^{(\alpha)}) \, \dd s , \qquad
   &J_t^{(6)} := \int_0^t
               \frac{b(s)^2}{\sigma(s)^2} \, r(Y_s^{(\alpha)} )^2 \, \dd s .
 \end{align*}
By \eqref{joint}, using Slutsky's lemma and the continuous mapping theorem,
 in order to prove the statement, it is sufficient to show
 \begin{gather}
  \lim_{t \uparrow T} \frac{I_{Y^{(\alpha)}}(t)}{I_{X^{(\alpha)}}(t)} = 1 ,
   \label{FI} \\[1mm]
   \text{$\frac{J_t^{(j)}}{\sqrt{I_{X^{(\alpha)}}(t)}} \stoch 0\quad$ \  as \ $t \uparrow T$
        \ for \ $j=1,2$,} \label{J12} \\[1mm]
  \text{$\frac{J_t^{(j)}}{I_{X^{(\alpha)}}(t)} \stoch 0\quad$ \ as \ $t \uparrow T$ \ for
        \ $j=3,4,5,6$.} \label{J36}
 \end{gather}
Using \eqref{Fisherp8} and the fact that \ $b(t)\ne0$ \ for all \ $t\in[t_0,T)$, \ we can apply
 L'Hospital's rule and we obtain,
 \[
   \lim_{t \uparrow T} \frac{I_{Y^{(\alpha)}}(t)}{I_{X^{(\alpha)}}(t)}
   = \lim_{t \uparrow T} \frac{\EE a(Y_t^{(\alpha)})^2}{\EE (X_t^{(\alpha)})^2} .
 \]
Using again the decomposition \eqref{DECO}, we have
 \[
   \EE a(Y_t^{(\alpha)})^2
   = \EE (X_t^{(\alpha)})^2 + 2 \EE X_t^{(\alpha)} \Delta_t^{(\alpha)}
     + \EE (\Delta_t^{(\alpha)})^2  + 2 \EE Y_t^{(\alpha)} r(Y_t^{(\alpha)})
     + \EE r(Y_t^{(\alpha)})^2 , \qquad t\in[0,T).
 \]
By Cauchy-Schwartz's inequality,
 \begin{align}
  &\frac{\big|\EE X_t^{(\alpha)} \Delta_t^{(\alpha)}\big|}{\EE (X_t^{(\alpha)})^2}
    \leq \sqrt{\frac{\EE (\Delta_t^{(\alpha)})^2}{\EE (X_t^{(\alpha)})^2}} ,
    \label{CS1} \\[2mm]
  &\frac{\big|\EE Y_t^{(\alpha)} r(Y_t^{(\alpha)})\big|}{\EE (X_t^{(\alpha)})^2}
   \leq \sqrt{\frac{2\big(\EE (X_t^{(\alpha)})^2 + \EE (\Delta_t^{(\alpha)})^2\big)}
                  {\EE (X_t^{(\alpha)})^2} \,
             \frac{\EE r(Y_t^{(\alpha)})^2}{\EE (X_t^{(\alpha)})^2}} , \label{CS2}
 \end{align}
 thus, in order to show \eqref{FI}, it is enough to check
 \begin{gather}
  \lim_{t \uparrow T} \frac{\EE (\Delta_t^{(\alpha)})^2}{\EE (X_t^{(\alpha)})^2} = 0 ,
   \label{Delta} \\
  \lim_{t \uparrow T} \frac{\EE r(Y_t^{(\alpha)})^2}{\EE (X_t^{(\alpha)})^2} = 0 .
   \label{rX}
 \end{gather}
In order to show \eqref{J12}, it is enough to prove
 \[
   \text{$\frac{J_t^{(j)}}{\sqrt{I_{X^{(\alpha)}}(t)}} \qmean 0\qquad$ \ as \ $t \uparrow T$
        \ for \ $j=1,2$,}
 \]
 which is equivalent to
 \begin{gather*} %\label{J1}
   \lim_{t \uparrow T}
    \frac{1}{I_{X^{(\alpha)}}(t)}
    \int_0^t \frac{b(s)^2}{\sigma(s)^2} \, \EE (\Delta_s^{(\alpha)})^2 \, \dd s
   = 0 ,  \\
   \lim_{t \uparrow T}
    \frac{1}{I_{X^{(\alpha)}}(t)}
    \int_0^t \frac{b(s)^2}{\sigma(s)^2} \, \EE r(Y_s^{(\alpha)})^2 \, \dd s
   = 0 . %\label{J2}
 \end{gather*}
By L'Hospital's rule,
 \begin{gather*}
  \lim_{t \uparrow T}
   \frac{1}{I_{X^{(\alpha)}}(t)}
   \int_0^t \frac{b(s)^2}{\sigma(s)^2} \, \EE (\Delta_s^{(\alpha)})^2 \, \dd s
  = \lim_{t \uparrow T} \frac{\EE (\Delta_t^{(\alpha)})^2}{\EE (X_t^{(\alpha)})^2} , \\
  \lim_{t \uparrow T}
   \frac{1}{I_{X^{(\alpha)}}(t)}
   \int_0^t \frac{b(s)^2}{\sigma(s)^2} \, \EE r(Y_s^{(\alpha)})^2 \, \dd s
  = \lim_{t \uparrow T} \frac{\EE r(Y_t^{(\alpha)})^2}{\EE (X_t^{(\alpha)})^2} ,
 \end{gather*}
 hence \eqref{J12} also follows from \eqref{Delta} and \eqref{rX}.
In order to show \eqref{J36}, it is enough to prove
 \[
   \text{$\frac{J_t^{(j)}}{I_{X^{(\alpha)}}(t)} \mean 0\qquad$ \ as \ $t \uparrow T$
         \ for \ $j=3,4,5,6$.}
 \]
For \ $j=4$ \ and \ $j=6$, \ by the previous argument, this follows directly from \eqref{Delta} and
 \eqref{rX}. For \ $j=3$ \ and \ $j=5$, \ this also follows from \eqref{Delta} and \eqref{rX}
 applying \eqref{CS1} and \eqref{CS2}.

The aim of the following discussions is to check \eqref{Delta} and \eqref{rX}.

First we consider the case \ $\alpha>0$ \ and \ $C>0$.
\ Then \ $b(t)>0$, \ $t\in[t_1,T)$, \ and, by Lemma \ref{LEMMA_EKVIVALENS_SINGULAR}, condition \eqref{Fisherp5}
 is equivalent to \ $\lim_{t\uparrow T}\int_0^tb(s)\,\dd s=\infty.$
\ Let us introduce the stochastic process
 \begin{align}\label{SEGED_PERTURB3}
   Z_t^{(\alpha)}
   := Y_t^{(\alpha)} \exp\left\{ - \alpha \int_0^t b(u) \, \dd u \right\} ,
   \qquad t \in [0,T) .
 \end{align}
Using \ $Y_t^{(\alpha)} = X_t^{(\alpha)} + \Delta_t^{(\alpha)}$, \ $t \in [0,T)$, \ and
 the equations \eqref{SolX} and \eqref{SolDelta}, we get
 \[
   Z_t^{(\alpha)}
   = \int_0^t
      \sigma(s) \exp\left\{ - \alpha \int_0^s b(u) \, \dd u \right\} \dd B_s
     + \alpha
       \int_0^t
        r(Y_s^{(\alpha)}) b(s) \exp\left\{ - \alpha \int_0^s b(u) \, \dd u \right\}
        \dd s
 \]
 for all \ $t \in [0,T)$.
\ Consequently, for all \ $t\in[0,T)$,
 \begin{equation} \label{Z2}
  \begin{split}
   (Z_t^{(\alpha)})^2
   &\leq 2\left(\int_0^t
                 \sigma(s)
                 \exp\left\{ - \alpha \int_0^s b(u) \, \dd u \right\}
                 \dd B_s \right)^2 \\
   &\phantom{\quad}
        + 2 \alpha^2
          \left(\int_0^t
                 r(Y_s^{(\alpha)}) b(s)
                 \exp\left\{ - \alpha \int_0^s b(u) \, \dd u \right\}
                 \dd s \right)^2 .
  \end{split}
 \end{equation}
Clearly, by Karatzas and Shreve \cite[Proposition 3.2.10]{KarShr}, for all \ $t\in[0,T)$,
 \[
   \EE \left(\int_0^t
              \sigma(s)
              \exp\left\{ - \alpha \int_0^s b(u) \, \dd u \right\}
              \dd B_s \right)^2
   = \int_0^t
      \sigma(s)^2
      \exp\left\{ - 2 \alpha \int_0^s b(u) \, \dd u \right\} \dd s ,
 \]
 and
 \begin{multline*}
  \EE \left(\int_0^t
             r(Y_s^{(\alpha)}) b(s)
             \exp\left\{ - \alpha \int_0^s b(u) \, \dd u \right\}
             \dd s \right)^2 \\
  = \int_0^t \int_0^t
     \EE \big[ r(Y_u^{(\alpha)}) r(Y_v^{(\alpha)}) \big] b(u) b(v)
      \exp\left\{ - \alpha \int_0^u b(w) \, \dd w
                  - \alpha \int_0^v b(w) \, \dd w \right\}
      \dd u \, \dd v .
 \end{multline*}
Using that \ $\big|\EE \big[ r(Y_u^{(\alpha)}) r(Y_v^{(\alpha)}) \big]\big|
 \leq \sqrt{\EE r(Y_u^{(\alpha)})^2 \, \EE r(Y_v^{(\alpha)})^2}$ \
 for all \ $u,v\in[0,T)$, \ we obtain
 \begin{align}\label{SEGED_PERTURB4}
   \begin{split}
    \EE&\left(\int_0^t
             r(Y_s^{(\alpha)}) b(s)
             \exp\left\{ - \alpha \int_0^s b(u) \, \dd u \right\}
             \dd s \right)^2 \\
    &\quad\leq \left(\int_0^t
              \sqrt{\EE r(Y_s^{(\alpha)})^2} \, |b(s)|
              \exp\left\{ - \alpha \int_0^s b(u) \, \dd u \right\}
              \dd s \right)^2,
     \qquad t\in[0,T).
   \end{split}
 \end{align}
Since \ $|b|^\gamma \leq 1 + |b|$ \ for all \ $b \in \RR$ \ and for all \ $\gamma\in[0,1)$,
 \ we have \ $|b|^{2\gamma} \leq 2(1 + b^2)$ \ for all \ $b \in \RR$ \ and for all
 \ $\gamma\in[0,1)$, \ and then, by \eqref{CONDITION_ON_r}, \ we obtain
 \begin{align*}
   r(Y_s^{(\alpha)})^2
   &\leq 2 L^2 (1 + |Y_s^{(\alpha)}|^{2\gamma})
    = 2L^2\left(1+\vert Z_s^{(\alpha)}\vert^{2\gamma}
                  \exp\left\{ 2 \gamma \alpha \int_0^s b(u) \, \dd u \right\}\right)\\
   &\leq 2 L^2
        \left(1 + 2\left( 1 + (Z_s^{(\alpha)})^2 \right)
                  \exp\left\{ 2 \gamma \alpha \int_0^s b(u) \, \dd u \right\}
        \right),\qquad s\in[0,T).
 \end{align*}
Recall that condition \eqref{Fisherp8} implied the existence of \ $c_1>0$, \ $c_2>0$ \ and
 \ $t_1\in[t_0,T)$ \ such that \eqref{liminfsup} holds.
By \eqref{Fisherp6}, there exists \ $t_2 \in [t_1,T)$ \ such
 that \ $\int_0^s b(u) \, \dd u > 0$ \ for all \ $s \in [t_2, T)$.
\ We check that
 \begin{align}\label{SEGED_PERTURB1}
   r(Y_s^{(\alpha)})^2
   \leq c\left(1 + (Z_s^{(\alpha)})^2 \right)
         \exp\left\{ 2 \gamma \alpha \int_0^s b(u) \, \dd u \right\} ,
         \qquad s \in [0,T) ,
 \end{align}
 with some appropriate \ $c>0$.
\ If \ $s\in[0,t_2]$, \ then
 \begin{align*}
    r(Y_s^{(\alpha)})^2
      &\leq 2 L^2
        \left(1 + 2\left( 1 + (Z_s^{(\alpha)})^2 \right)
                  \exp\left\{ 2 \gamma \alpha \sup_{v\in[0,t_2]}\int_0^v \vert b(u)\vert\,\dd u\right\}
        \right)\\[1mm]
     &\leq 8L^2\left( 1 + (Z_s^{(\alpha)})^2 \right)
         \exp\left\{ 2 \gamma \alpha \sup_{v\in[0,t_2]}\int_0^v \vert b(u)\vert\,\dd u\right\}\\[1mm]
     &\leq c'\left( 1 + (Z_s^{(\alpha)})^2 \right)
          \exp\left\{ 2 \gamma \alpha \int_0^s b(u) \, \dd u \right\},
 \end{align*}
 where
 \[
    c':=\frac{8L^2\exp\left\{ 2 \gamma \alpha \sup_{v\in[0,t_2]}\int_0^v \vert b(u)\vert\,\dd u\right\}}
          {\exp\left\{ 2 \gamma \alpha \inf_{v\in[0,t_2]}\int_0^v b(u)\,\dd u\right\}} .
 \]
If \ $s\in(t_2,T)$, \ then
 \[
     r(Y_s^{(\alpha)})^2
       \leq 8L^2\left( 1 + (Z_s^{(\alpha)})^2 \right)
            \exp\left\{ 2 \gamma \alpha \int_0^s b(u) \, \dd u \right\}.
 \]
Hence \eqref{SEGED_PERTURB1} is satisfied with \ $c:=\max\{c',8L^2\}$.
\ Thus, using that \ $\sqrt{a+b}\leq\sqrt{a}+\sqrt{b}$ \ for all \ $a,b\geq0$, \ we have
 \begin{equation} \label{r}
  \sqrt{\EE r(Y_s^{(\alpha)})^2}
  \leq \sqrt{c}\left(1 + \sqrt{\EE(Z_s^{(\alpha)})^2} \right)
         \exp\left\{ \gamma \alpha \int_0^s b(u) \, \dd u \right\} , \qquad
  s \in [0,T).
 \end{equation}
Applying \eqref{Z2} and \eqref{r}, we obtain for all \ $t\in[0,T)$,
 \begin{align*}
  \EE (Z_t^{(\alpha)})^2
  &\leq 2 \int_0^t
          \sigma(s)^2
          \exp\left\{ - 2 \alpha \int_0^s b(u) \, \dd u \right\} \dd s \\
  &\phantom{\quad\,}
   + 2c\alpha^2
     \left( \int_0^t \left(1 + \sqrt{\EE(Z_s^{(\alpha)})^2} \right) |b(s)|
             \exp\left\{ - (1-\gamma) \alpha \int_0^s b(u) \, \dd u \right\}
             \dd s \right)^2 .
 \end{align*}
By the choice of \ $t_1$, \ we have for all \ $t\in[t_1,T)$,
 \begin{align*}
  \int_0^t |b(s)|
           \exp\left\{ - (1-\gamma) \alpha \int_0^s b(u) \, \dd u \right\}
           \dd s
  &\leq \int_0^{t_1} |b(s)|
         \exp\left\{ - (1-\gamma) \alpha \int_0^s b(u) \, \dd u \right\}
         \dd s \\
  & \phantom{\quad}
        + \int_{t_1}^t b(s)
           \exp\left\{ - (1-\gamma) \alpha \int_0^s b(u) \, \dd u \right\}
           \dd s ,
 \end{align*}
 and for all \ $t\in[t_1,T)$, \
 \[
   \int_{t_1}^t b(s)
    \exp\left\{ - (1-\gamma) \alpha \int_0^s b(u) \, \dd u \right\} \dd s
   \leq \frac{1}{(1-\gamma) \alpha}
        \exp\left\{ - (1-\gamma) \alpha \int_0^{t_1} b(u) \, \dd u \right\}.
 \]
Hence there exists \ $K_1 > 0$ \ such that
 \begin{equation} \label{K1}
   \int_0^t |b(s)|
    \exp\left\{ - (1-\gamma) \alpha \int_0^s b(u) \, \dd u \right\} \dd s
   \leq K_1 , \qquad t \in [0,T) .
 \end{equation}
Using \eqref{liminfsup}, \eqref{K1} and that
\ $\sqrt{a+b}\leq\sqrt{a}+\sqrt{b}$ \ for all \ $a,b\geq 0$, \ we obtain
 \[
   \sqrt{\EE (Z_t^{(\alpha)})^2}
   \leq \psi_1(t) + \int_{t_2}^t \psi_2(u) \sqrt{\EE (Z_u^{(\alpha)})^2} \, \dd u ,
   \qquad t \in [t_2,T) ,
 \]
 where
 \begin{align*}
  \psi_1(t)
  &:= \left(2\int_0^{t_2}\sigma(s)^2\exp\left\{-2\alpha\int_0^sb(u)\,\dd u\right\}\,\dd s\right)^{1/2}
      +\left(2 c_2 \int_{t_2}^t b(s) \, \dd s\right)^{1/2}
      +\sqrt{2c}\alpha K_1\\
  &\phantom{:=\;}
      +\int_0^{t_2} \psi_2(u) \sqrt{\EE (Z_u^{(\alpha)})^2}\,\dd u\\[1mm]
  &=:\left(2 c_2 \int_{t_2}^t b(s) \, \dd s\right)^{1/2}+K_2,
     \qquad t\in[t_2,T),\\[2mm]
  \psi_2(t)
  &:= \sqrt{2c}\alpha |b(t)|
      \exp\left\{ -(1-\gamma)\alpha\int_0^t b(u) \, \dd u \right\},
      \qquad t\in[t_2,T).
 \end{align*}
Hence, by Lemma \ref{Gronwall} (generalized Gr\"onwall's inequality), we obtain
 \begin{align*}
  \sqrt{\EE (Z_t^{(\alpha)})^2}
  &\leq K_2\exp\left\{\int_{t_2}^t \psi_2(u)\,\dd u\right\}
     +\int_{t_2}^t\psi_1'(u)\exp\left\{\int_u^t \psi_2(v)\,\dd v\right\}\,\dd u,
          \qquad t\in[t_2,T).
 \end{align*}
Using that, by \eqref{K1},
 \[
  \exp\left\{\int_{t_2}^t \psi_2(u)\,\dd u\right\}
      \leq \ee^{\sqrt{2c}\alpha K_1}=:K_3,
      \qquad t\in[t_2,T),
 \]
 and that
 \[
   \psi_1'(t)
      =c_2b(t)\left(2 c_2 \int_{t_2}^t b(s) \, \dd s\right)^{-1/2}>0,
      \qquad t\in[t_2,T),
 \]
 we have
 \begin{align*}
   \sqrt{\EE (Z_t^{(\alpha)})^2}
    &\leq K_2K_3+K_3\int_{t_2}^t\psi_1'(u)\,\dd u
     = K_2K_3+K_3\left(2 c_2 \int_{t_2}^t b(s) \, \dd s\right)^{1/2},
    \qquad t\in[t_2,T).
 \end{align*}
Then, using \eqref{Fisherp6}, we have there exist \ $c''>0$ \ and \ $t_3\in(t_2,T)$
 \ such that
 \[
   \sqrt{\EE (Z_t^{(\alpha)})^2}
   \leq c'' \left(\int_{t_2}^t b(s) \, \dd s\right)^{1/2} , \qquad
     t\in[t_3, T).
 \]
By Lyapunov's inequality, since \ $0\leq 2\gamma<2$,
 \begin{align*} %\label{Z}
  \EE |Z_t^{(\alpha)}|^{2\gamma}
  \leq \left(\EE (Z_t^{(\alpha)})^2\right)^\gamma
  \leq (c'')^{2\gamma} \left(\int_{t_2}^t b(s) \, \dd s \right)^\gamma , \qquad
  t \in [t_3,T),
 \end{align*}
 and thus, by \eqref{SEGED_PERTURB3},
 \begin{align} \label{Y}
  \begin{split}
   \EE |Y_t^{(\alpha)}|^{2\gamma}
     &=\exp\left\{2\gamma\alpha\int_0^t b(u)\,\dd u \right\}
        \EE |Z_t^{(\alpha)}|^{2\gamma}\\
     &\leq (c'')^{2\gamma} \exp\left\{ 2 \gamma \alpha \int_0^t b(u) \, \dd u \right\}
         \left(\int_{t_2}^t b(s) \, \dd s\right)^\gamma , \qquad t \in [t_3, T).
   \end{split}
 \end{align}
Applying again \eqref{liminfsup}, for all \ $t\in[t_1,T)$, \ we have
\begin{align}\label{SEGED_PERTURB5}
 \begin{split}
   &\EE(X_t^{(\alpha)})^2
       =\exp\left\{2\alpha \int_0^t b(v)\,\dd v\right\}
          \int_0^t\sigma(u)^2\exp\left\{-2\alpha \int_0^u b(v)\,\dd v\right\}\,\dd u\\
   &\quad \geq \exp\left\{2\alpha \int_0^t b(v)\,\dd v\right\}
            \left(\int_0^{t_1}\sigma(u)^2\exp\left\{-2\alpha \int_0^u b(v)\,\dd v\right\}\,\dd u
                   +\int_{t_1}^tc_1b(u)\,\dd u\right).
 \end{split}
\end{align}
Hence, using also \eqref{CONDITION_ON_r}, \ for all \ $t\in[t_3,T)$, \ we have
 \begin{align*}
  \frac{\EE r(Y_t^{(\alpha)})^2}{\EE (X_t^{(\alpha)})^2}
   & \leq \frac{2 L^2 (1+\EE |Y_t^{(\alpha)}|^{2\gamma})}{\EE (X_t^{(\alpha)})^2} \\
   & \leq \frac{2 L^2
              \Big(1+(c'')^{2\gamma} \exp\left\{ 2 \gamma \alpha
                                    \int_0^t b(u) \, \dd u \right\}
                      \left(\int_{t_2}^t b(s) \, \dd s\right)^\gamma \Big)}
             {\exp\left\{ 2 \alpha \int_0^t b(u) \, \dd u \right\}
                \left(\int_0^{t_1}\sigma(u)^2\exp\left\{-2\alpha\int_0^u b(v)\,\dd v\right\}\,\dd u
                        + c_1\int_{t_1}^t b(s)\,\dd s\right)}.
 \end{align*}
Using that \ $\gamma\in[0,1)$, \ $c_1>0$, \ and \ $\lim_{x\to\infty}\frac{x^a}{\ee^{bx}}=0$
 \ for all \ $a>0$ \ and \ $b>0$, \ by \eqref{Fisherp6}, we get \eqref{rX}.

\noindent Now we turn to prove \eqref{Delta}.
Using \eqref{SolDelta}, by the same way that we derived \eqref{SEGED_PERTURB4}, one can get
 \begin{align*}
   \EE(\Delta_t^{(\alpha)})^2
      \leq \alpha^2
         \left(\int_0^t \sqrt{\EE r(Y_s^{(\alpha)})^2} \, |b(s)|
               \exp\left\{\alpha\int_s^t b(u)\,\dd u\right\}\dd s\right)^2,
               \qquad t\in[0,T).
 \end{align*}
By \eqref{Fisherp6}, there exists \ $t_4\in(t_3,T)$ \ such that
 \[
   \exp\left\{ 2 \gamma \alpha \int_0^t b(u) \, \dd u \right\}
                      \left(\int_{t_2}^t b(s) \, \dd s\right)^\gamma
   \geq 1 ,
                      \qquad t\in[t_4,T),
 \]
Hence, using \eqref{CONDITION_ON_r} and \eqref{Y}, we have
 \begin{align*}
    \sqrt{\EE r(Y_t^{(\alpha)})^2}
       \leq \sqrt{2}L\sqrt{1+(c'')^{2\gamma}}
             \exp\left\{\gamma\alpha\int_0^t b(u)\,\dd u\right\}
             \left(\int_{t_2}^t b(u)\,\dd u\right)^{\gamma/2},
        \qquad t \in [t_4, T),
 \end{align*}
 and then for all \ $t\in[t_4,T)$,
 \begin{align*}
   \EE(\Delta_t^{(\alpha)})^2
     &\leq 2\alpha^2
          \left(\int_0^{t_4} \!\!\!\! \sqrt{\EE r(Y_s^{(\alpha)})^2} \, |b(s)|
               \exp\left\{\alpha \!\! \int_s^t b(u)\,\dd u\right\}\dd s\right)^2\\
     &%\phantom{\leq}
           +2\alpha^22L^2(1+(c'')^{2\gamma})
              \!\!\left(\int_{t_4}^t \!\! \vert b(s)\vert
                     \exp\left\{\!\gamma\alpha \!\! \int_0^s \!\! b(u)\,\dd u
                      +\alpha \!\! \int_s^t \!\! b(u)\,\dd u\right\}
                     \!\!\left(\int_{t_2}^s \!\! b(u)\,\dd u\right)^{\gamma/2}\!\dd s
             \right)^2.
 \end{align*}
Hence, by \eqref{SEGED_PERTURB5}, we have
 \begin{align*}
  \frac{\sqrt{\EE (\Delta_t^{(\alpha)})^2}}{\sqrt{\EE (X_t^{(\alpha)})^2}}
    &\leq \lim_{t \uparrow T}
       \frac{\sqrt{2}\alpha\int_0^{t_4}\sqrt{\EE r(Y_s^{(\alpha)})^2}\vert b(s)\vert
             \exp\left\{\alpha\int_s^t b(u)\,\dd u\right\}\,\dd s}
      {\exp\left\{\alpha\int_0^t b(u)\,\dd u\right\}
          \left(\int_0^{t_1}\sigma(u)^2\exp\left\{-2\alpha\int_0^u b(v)\,\dd v\right\}\,\dd u
               +c_1\int_{t_1}^t b(s)\,\dd s\right)^{\frac{1}{2}}}\\[2mm]
  &\phantom{\leq\;}
   +\frac{2\alpha L\sqrt{1+(c'')^{2\gamma}}
           \int_{t_4}^t\vert b(s)\vert
                    \exp\left\{-(1-\gamma)\alpha\int_0^s b(u)\,\dd u\right\}
                     \left(\int_{t_2}^s b(u)\,\dd u\right)^{\frac{\gamma}{2}}\,\dd s}
          {\left(\int_0^{t_1}\sigma(u)^2\exp\left\{-2\alpha\int_0^u b(v)\,\dd v\right\}\,\dd u
                        + c_1\int_{t_1}^t b(s)\,\dd s\right)^{1/2}},
 \end{align*}
and then, by L'Hospital's rule, we conclude
 \begin{align*}
  \lim_{t \uparrow T}\frac{\sqrt{\EE (\Delta_t^{(\alpha)})^2}}{\sqrt{\EE (X_t^{(\alpha)})^2}}
  &\leq\frac{2\alpha L\sqrt{1+(c'')^{2\gamma}}}{\sqrt{c_1}}
     \lim_{t \uparrow T}
       \frac{\int_{t_4}^t\vert b(s)\vert
                    \exp\left\{-(1-\gamma)\alpha\int_0^s b(u)\,\dd u\right\}
                     \left(\int_{t_2}^s b(u)\,\dd u\right)^{\frac{\gamma}{2}}\,\dd s}
            {\left(\int_{t_1}^t b(s)\,\dd s\right)^{1/2}}\\
  &=\frac{2\alpha L\sqrt{1+(c'')^{2\gamma}}}{\sqrt{c_1}}
     \lim_{t \uparrow T}
      \frac{\vert b(t)\vert
             \exp\left\{-(1-\gamma)\alpha\int_0^t b(u)\,\dd u\right\}
                     \left(\int_{t_2}^t b(u)\,\dd u\right)^{\frac{\gamma}{2}}}
            {\frac{1}{2}\left(\int_{t_1}^t b(s)\,\dd s\right)^{-1/2}b(t)}
     =0,
 \end{align*}
 since \ $\lim_{t \uparrow T}\int_0^t b(s)\,\dd s=\infty$ \ and \ $\lim_{x\to\infty}\frac{x^a}{\ee^{bx}}=0$
 \ for all \ $a>0$ \ and \ $b>0$.

Now we consider case \ $\alpha = 0$.
\ Then the SDE \eqref{perturbed_SDE} and the SDE \eqref{special_SDE} have the same unique
 strong solution
 \begin{align*} % \label{SolY}
  Y^{(0)}_t = \int_0^t \sigma(s) \, \dd B_s , \qquad t \in [0,T) .
 \end{align*}
Hence \ $\Delta_t^{(0)}=Y_t^{(0)}-X_t^{(0)}=0$, \ $t\in[0,T)$, \ and then we get
 \[
   \halpha_t^{(Y^{(0)})}
   = \frac{\int_0^t \frac{b(s)}{\sigma(s)} \, X_s^{(0)} \, \dd B_s + J_t^{(2)}}
          {\int_0^t \frac{b(s)^2}{\sigma(s)^2} \, ( X_s^{(0)} )^2 \, \dd s
           + J_t^{(5)} + J_t^{(6)}} ,
 \]
 holds asymptotically as \ $t\uparrow T$ \ with probability one.
Note that \eqref{Delta} is satisfied, since \ $\Delta_t^{(0)}=0$, \ $t\in[0,T)$.
\ Hence, in order to prove the statement, it is enough to check \eqref{rX}.
Clearly,
 \[
    X_t^{(0)}
   \distre Y_t^{(0)}
   \distre \cN\left( 0, \int_0^t \sigma(s)^2 \, \dd s \right) , \qquad
   t \in [0,T) ,
 \]
 and hence
 \begin{align*}
  &\EE (X_t^{(0)})^2
     = \int_0^t\sigma(s)^2\,\dd s,\qquad t\in[0,T),\\
  &\EE |Y_t^{(0)}|^{2 \gamma}
     = \left( \int_0^t \sigma(s)^2 \, \dd s \right)^\gamma
        \EE\vert\xi\vert^{2\gamma},
      \qquad t\in[0,T),
 \end{align*}
 where \ $\xi$ \ is a standard normally distributed random variable.
Then, by \eqref{CONDITION_ON_r}, \ we have
 \begin{align*}
   \lim_{t\uparrow T}\frac{\EE r(Y^{(0)}_t)^2}{\EE (X^{(0)}_t)^2}
        \leq \lim_{t\uparrow T}\frac{2L^2(1+\EE\vert Y^{(0)}_t\vert^{2\gamma})}{\EE (X^{(0)}_t)^2}
         = \lim_{t\uparrow T}\frac{2L^2\left(1+\EE\vert\xi\vert^{2\gamma}
                         \left(\int_0^t\sigma(s)^2\,\dd s\right)^\gamma\right)}
                               {\int_0^t\sigma(s)^2\,\dd s}
         =0,
 \end{align*}
 where the last step can be checked as follows.
By \eqref{Fisherp8}, \ $\lim_{t\uparrow T}\frac{\vert b(t)\vert}{\sigma(t)^2}=\vert C\vert\in(0,\infty)$,
 \ and hence there exist \ $c_2>0$ \ and \ $t_1\in[t_0,T)$ \ such that
 \ $\vert b(t)\vert <c_2\sigma(t)^2$ \ for all \ $t\in[t_1,T)$.
\ Then \eqref{Fisherp6} yields \ $\lim_{t\uparrow T}\int_0^t\sigma(s)^2\,\dd s=\infty$ \
 concluding the proof of the present case.

Finally, we consider the case \ $\alpha<0$ \ and \ $C<0$.
\ For all \ $\beta\in\RR$, \ let us consider the process \ $(V_t^{(\beta)})_{t\in[0,T)}$ \ given
 by the SDE
 \begin{align*}
  \begin{cases}
   \dd V_t^{(\beta)}
   = \beta\,\widetilde b(t)a(V_t^{(\beta)}) \, \dd t
     + \sigma(t) \, \dd B_t,\qquad t\in[0,T),\\
   \phantom{\dd} V_0^{(\beta)}=0,
  \end{cases}
 \end{align*}
 where \ $\widetilde b(t):=-b(t)$, \ $t\in[0,T)$.
\ Note that if conditions \eqref{Fisherp5} and \eqref{Fisherp8} are satisfied with functions
 \ $b$ \ and \ $\sigma$ \ and with parameters \ $\alpha<0$ \ and \ $C<0$, \ then they are also
 satisfied with the functions \ $\widetilde b=-b$ \ and \ $\sigma$ \ and with parameters
 \ $-\alpha>0$ \ and \ $-C>0$.
Hence
 \[
   \sqrt{I_{V^{(-\alpha)}}(t)} \, (\widehat{(-\alpha)}_t^{(V^{(-\alpha)})} - (-\alpha) )
   \distr \frac{1}{\sqrt{2}} \,
          \frac{\int_0^1 W_s \, \dd W_s}{\int_0^1 (W_s)^2 \, \dd s} \qquad
   \text{as \ $t \uparrow T$.}
 \]
By the uniqueness of a strong solution, the process \ $(Y_t^{(\alpha)})_{t\in[0,T)}$
 \ given by the SDE \eqref{perturbed_SDE} and the process \ $(V_t^{(-\alpha)})_{t\in[0,T)}$
 \ coincide, and hence
 \ $I_{V^{(-\alpha)}}(t)=I_{Y^{(\alpha)}}(t)$ \ for all \ $t\in(0,T)$,
 \ and \ $\widehat{(-\alpha)}_t^{(V^{(-\alpha)})}=-\halpha_t^{(Y^{(\alpha)})}$
 \ holds asymptotically as \ $t\uparrow T$ \ with probability one concluding the proof.
\proofend

\begin{Rem}
We note that in the proof of Theorem \ref{sing_perturb} instead of Lemma \ref{Gronwall}
 (a generalization of Gr\"onwall's inequality) we could use Bainov and Simeonov
 \cite[Theorem 1.3]{BaiSem}, which is another generalization of Gr\"onwall's inequality.
 But the calculations would be more complicated without any improvement or refinement of the result.
\end{Rem}

\end{document}